\documentclass[a4paper]{elsarticle}

\usepackage{enumerate}
\usepackage{hhline}
\usepackage{tabulary}
\newcolumntype{K}[1]{>{\centering\arraybackslash}p{#1}}
\newcolumntype{M}[1]{>{\centering\arraybackslash}m{#1}}
\usepackage{caption}
\usepackage[labelformat=simple]{subcaption}
\usepackage{multirow}
\usepackage{array}

\def\eop{\vrule width4pt height5pt depth6pt}

\journal{arXiv}

\usepackage{lineno,hyperref}
\usepackage{amssymb,amsmath,amscd,mathrsfs,bm,accents,amsthm,graphicx,psfrag,epstopdf}
\newtheorem{theorem}{\sc Theorem}[]

\newtheorem{lemma}[theorem]{\sc Lemma}

\theoremstyle{definition}
\theoremstyle{definition}
\theoremstyle{definition}
\theoremstyle{definition}
\theoremstyle{definition}

\def\bit{\begin{itemize}}   \def\eit{\end{itemize}}
\def\ben{\begin{enumerate}} \def\een{\end{enumerate}}
\def\bpic{\begin{picture}}  \def\epic{\end{picture}}
\def\bfig{\begin{figure}}   \def\efig{\end{figure}}
\def\bc{\begin{center}}     \def\ec{\end{center}}
\def\beq{\begin{equation}}  \def\eeq{\end{equation}}
\def\beqn{\begin{eqnarray}} \def\eeqn{\end{eqnarray}}
\def\beqnn{\begin{eqnarray*}}   \def\eeqnn{\end{eqnarray*}}
\def\barr{\begin{array}}    \def\earr{\end{array}}
\def\tb{\textbf}
\def\dt{\Delta(T)}

\begin{document}

\begin{frontmatter}



\title{The game chromatic index of trees of maximum degree 4 with at most three degree-four vertices in a row}


\author[add1]{Wai Lam Fong\corref{cor1}}
\ead{s1118833@s.eduhk.hk}
\author[add1]{Wai Hong Chan}
\ead{waihchan@eduhk.hk}
\cortext[cor1]{Corresponding author}
\address[add1]{Department of Mathematics and Information Technology, The Education University of Hong Kong, Tai Po, Hong Kong SAR, China.}

\begin{abstract}

Fong et al. (The game chromatic index of some trees with maximum degree four and adjacent degree-four vertices, J. Comb Optim 36 (2018) 1-12) proved that the game chromatic index of any tree $T$ of maximum degree 4 whose degree-four vertices induce a forest of paths of length $l$ less than 2 is at most 5. In this paper, we show that the bound 5 is also valid for $l\leq 2$. This partially solves the problem of characterization of the trees whose game chromatic index exceeds the maximum degree by at most 1, which was proposed by Cai and Zhu (Game chromatic index of $k$-degenerate graphs, J. Graph Theory 36 (2001) 144-155). 
\end{abstract}



\end{frontmatter}

\section{Introduction} We consider in this paper an edge-coloring game studied in~\cite{A2003, A2006,CN,EFHK,FC,FCN}. In the game, two players, Alice and Bob, alternately select a color from a set of colors and put it on an uncolored edge of an initially uncolored, finite and simple graph $G$ such that adjacent edges receive distinct colors. Alice wins the game if all edges of $G$ are colored finally; otherwise, Alice loses. Moreover, Bob begins and he is allowed to skip any number of turns throughout the game, while skipping is forbidden for Alice. The game played with a set of $k$ colors is called the {\it $k$-edge-coloring game}. If Alice has a winning strategy for the $k$-edge-coloring game played on a graph $G$, then $G$ is called {\it $k$-edge-game-colorable}. The parameter {\it game chromatic index} $\chi'_g(G)$ of a graph $G$, which was introduced by Cai and Zhu~\cite{CZ} for a related game, is defined as the least $k$ such that Alice has a winning strategy for the $k$-edge-coloring game played on $G$. Bodlaender~\cite{B} introduced a similar kind of coloring games, in which vertices are colored instead of edges, and the corresponding parameter is called the {\it game chromatic number}.

Cai and Zhu~\cite{CZ} proved that  $\chi'_g(T)\leq\dt +2$ for any tree $T$ of maximum degree $\dt$. Erd{\"o}s et al.~\cite{EFHK} showed that best possible upper bound for the class of trees $T$ with maximum degree $\dt$ is at least $\dt +1$ if $\dt\geq 2$. The upper bound $\dt +1$ holds when $\dt =3$~\cite{A2003, CZ} or $\dt\geq 5$~\cite{A2006, EFHK}. Any tree of maximum degree 2 is a path, which obviously means that all edges can be successfully colored in the 3-edge-coloring game played on any tree $T$ of $\dt =2$. As a result, $\dt =4$ is the only case that was left open. On the one hand, researchers have not proved that all trees of maximum degree 4 are 5-edge-game-colorable; on the other hand, no trees of maximum degree 4 that are not 5-edge-game-colorable have been found. Some researchers have been working on finding subclasses of trees $T$ of maximum degree 4 which are 5-edge-game-colorable, by considering the connectivity and distribution of 4-vertices (degree-four vertices) in $T$. For example, the subclass that the 4-vertices in $T$ induce a {\it linear forest}, which is a disjoint union of paths, were studied \citep{CN,FC,FCN}. These studies will be introduced below.

For trees $T$ with $\dt=4$. Chan and Nong~\cite{CN} proved that the upper bound $\dt +1$, i.e., 5, is also sharp when the set of 4-vertices of $T$ is independent, or when $T$ is a caterpillar, which may contain a long path of 4-vertices. After that, Fong and Chan~\cite{FC} found that any 4-edge-game-colorable tree with maximum degree 4, in which any 4-vertex is adjacent to at most two vertices of degree not less than 3, is 5-edge-game-colorable. The result of Fong and Chan~\cite{FC} not only provided a new 5-edge-game-colorable subclass of trees with maximum degree 4, but also confirmed that any tree is $k$-edge-game-colorable when $k$ is greater than its game chromatic index. This result also partially answered a basic yet challenging and open question raised by Zhu~\cite{Z}: Is having more colors always an advantage for Alice? Furthermore, Fong et al.~\cite{FCN} showed that the bound 5 is still valid when the subgraph of $T$ induced by all its 4-vertices is a forest of paths of length $l$ at most 1. Moreover, Fong et al.~\cite{FCN} conjectured that the result is true for any natural number $l$. In this paper, we confirm this conjecture for the case that $l\leq 2$. 

The game studied in this paper is under the game variant that Bob begins and has the right to skip any number of turns. This variant was proposed by Andres~\cite{A2006}. Andres~\cite{A2006} proposed six game variants in which either Alice or Bob could take the first move; and either Alice, Bob, or none of them could be allowed to skip. In other words, there are two and three choices of the player who begins and the player who can skip, respectively, so that there are six variants in total. Andres~\cite{A2006} showed that the game chromatic index under the variant studied in this paper is an upper bound on the index under the other five variants. This means if one wants to establish an upper bound on the game chromatic index of a class of graphs, then one may consider the variant that Bob is the player who begins and is allowed to skip. Andres et al.~\cite{AHS} studied the relationship between graph structure and the effect of the game variants on the game chromatic index of graphs.

The rest of this paper is structured as follows. In Section 2, we will present our terminology and notation, and the notion of decomposition of a tree into independent subtrees, which facilitate our proof of the conjecture for $l\leq 2$, i.e., our main theorem. Moreover, we will state Lemma~\ref{lem:inductive} and use it to prove the main theorem. At the end of Section 2, we will establish Lemma~\ref{lem:1permitted}, which will be employed to prove Lemma~\ref{lem:inductive} in Sections 3 and 4. Future work will be discussed in Section 5.

\section{Decomposition of a tree into independent subtrees}

\begin{theorem}\label{mainthm}
Let $T$ be a finite tree with $\dt=4$. If the subgraph of $T$ induced by all its 4-vertices is a forest of paths of length at most 2, then $\chi'_g(T) \leq 5$.
\end{theorem}

We now define the following terms for any colored or uncolored tree:

\begin{itemize}
\item A {\it k-vertex} is a vertex of degree $k$.
\item The edge incident with a leaf is called a {\it leaf-edge}.
\item The non-pendant vertex incident with a leaf-edge is called the {\it root of a leaf-edge}.
\item A {\it trivial path} is a path of length zero, i.e., a path consisting of one vertex.
\item A {\it star-node} is a vertex connected to at least three roots of colored leaf-edges by edge-disjoint (maybe trivial) paths.
\item Any star-node connected to precisely $k$ roots of colored leaf-edges by edge-disjoint (maybe trivial) paths is called a {\it k-SN}.
\item A {\it star-edge} is an edge incident with a star-node.
\item The path connecting two star-nodes is called a {\it star-path}.
\item A {\it k-leaf-edge-colored tree} ($k$-$LCT$)  is a tree containing precisely $k$ colored leaf-edges.
\item In a tree containing vertex $u$, a {\it u-branch} is a maximal subtree containing precisely one edge incident with $u$. A $u$-branch with no colored edges is called an {\it uncolored $u$-branch}.
\item A {\it leaf-path} is the path connecting the root of a colored leaf-edge with its nearest star-node.
\item {\it $S_3^4$} is the star with exactly three leaves, i.e., the claw, formed by four 4-vertices. 
\item {\it $P_n^4$} is the path of length $n-1$ formed by $n$ vertices of degree 4.
\end{itemize}

Note that $S_3^4$ and $P_4^4$ must not exist in a tree $T$ with $\dt=4$ and its subgraph induced by all its 4-vertices being a forest of paths of length at most 2. Also, we have the following remarks for any tree $T$ with $\dt=4$, based on the above definitions:
\begin{enumerate}
\item Any star-node is of degree 3 or 4.
\item A tree has zero star-nodes if and only if it has less than three colored leaf-edges.  A 3-LCT has precisely one 3-SN; a 4-LCT has either exactly two 3-SNs or exactly one 4-SN.
\end{enumerate}

\begin{figure}
\centering
\subcaptionbox{An example of Type 3 subtrees. The path $e_1e_2e_3$ is the subgraph induced by all the uncolored edges $e_1$, $e_2$ and $e_3$. Colors 4 and 5 are available for $e_1$ and $e_2$; and colors 3, 4 and 5 are available for $e_3$. \label{fig:3ex}} [.8\linewidth]
 {\includegraphics[scale=.5]{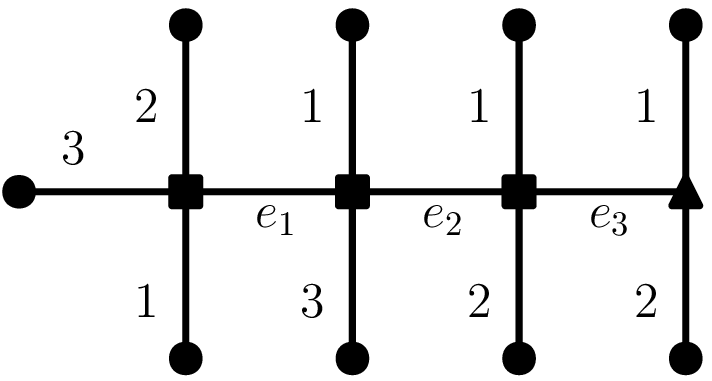}}
\subcaptionbox{An example of Type 4 subtrees. Colors 4 and 5 are available for $v_0v_1$ and $v_1v_2$; colors 3 and 5 are available for $v_2v_3$; and colors 1, 2, 3 and 5 are available for $v_3v_4$. \label{fig:4ex}} [.8\linewidth]
 {\includegraphics[scale=.5]{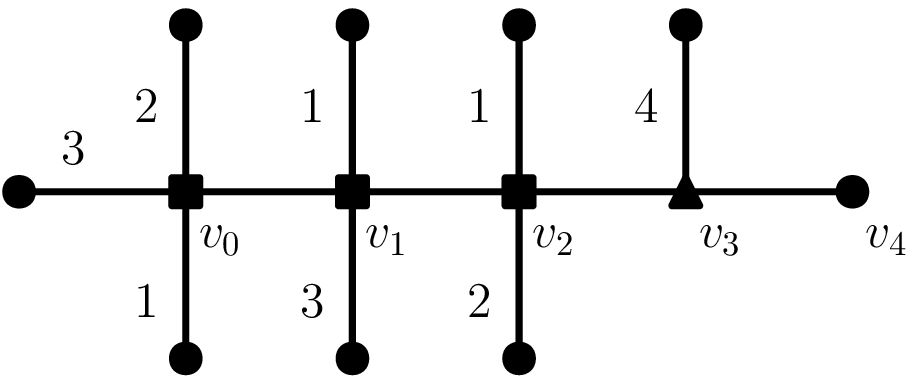}}
\subcaptionbox{Type 5. It is under the constraint that $d(w)=2$ and the trivial constraint that the two colored edges incident with $u$ do not have the same color.\label{fig:5}} [.3\linewidth]
 {\includegraphics[scale=.5]{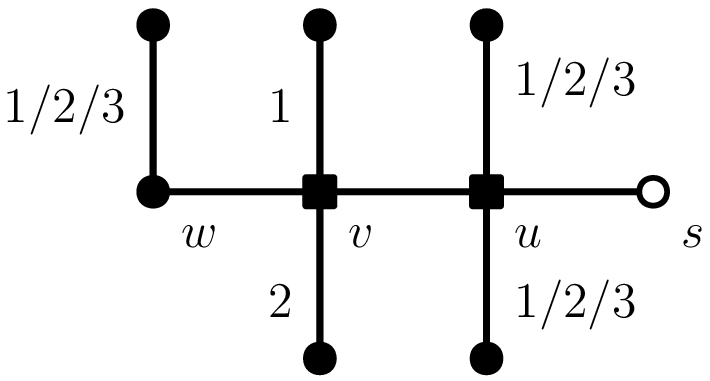}}\hfill%
\subcaptionbox{Type 6 \label{fig:6}} [.3\linewidth]
 {\includegraphics[scale=.5]{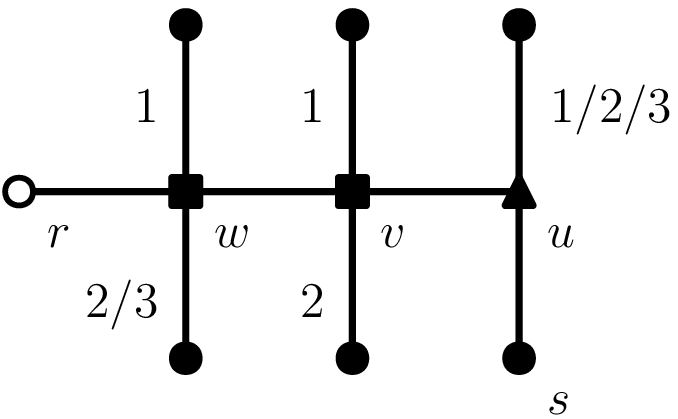}}\hfill%
\subcaptionbox{Type 7. Vertex $x$ exists if $d(r)=3$. \label{fig:7}} [.3\linewidth]
 {\includegraphics[scale=.5]{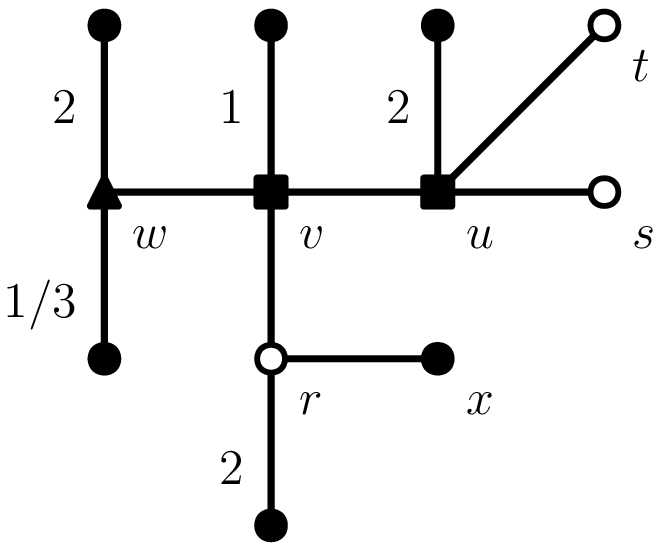}}
\subcaptionbox{Type 8 \label{fig:8}} [.3\linewidth]
 {\includegraphics[scale=.5]{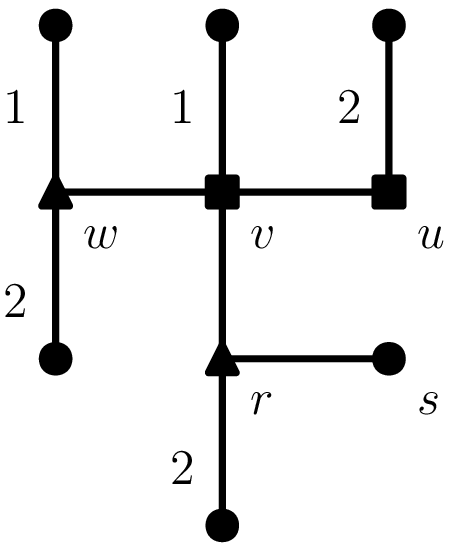}}\hfill%
\subcaptionbox{Type 9\label{fig:9}} [.3\linewidth]
 {\includegraphics[scale=.5]{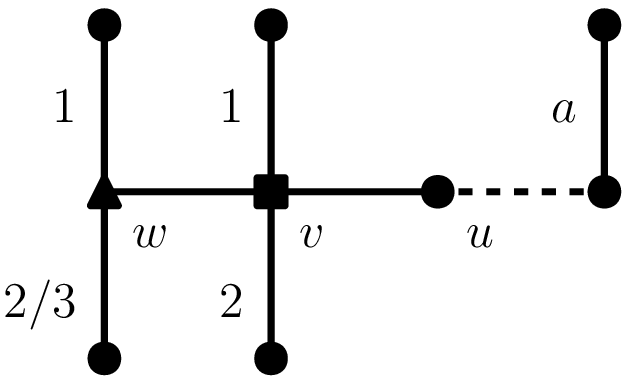}} \hfill%
\subcaptionbox{Type 10\label{fig:10}} [.3\linewidth]
 {\includegraphics[scale=.5]{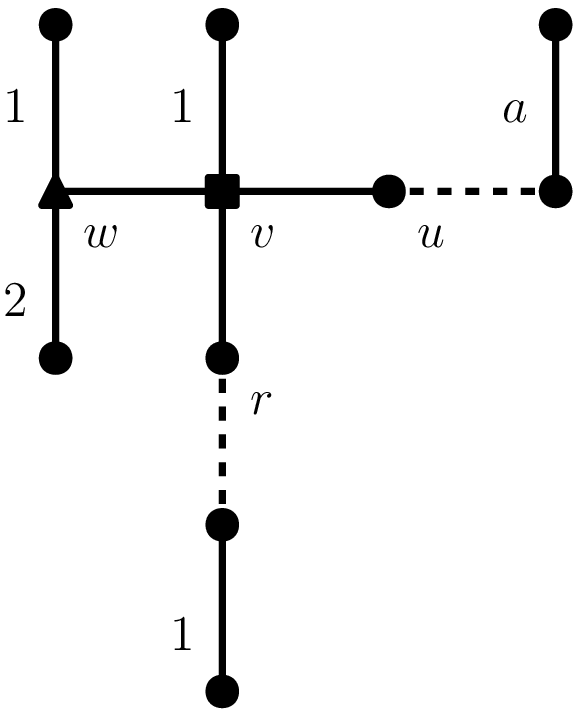}} 
 	\caption{Types 3 - 10.\label{fig:4-10}}
\end{figure}

The following are some notations and remarks about our figures. Let $\{1,2,3,4,5\}$ be the color set of the game with five colors, and letter {\it a, b, c} be arbitrary colors in the color set. A {\it filled rectangle}, a {\it filled triangle}, an {\it unfilled circle}, an {\it unfilled rhombus}, and a {\it filled circle} represent a 4-vertex, a 3-vertex, a vertex of degree at most 3, a vertex of degree at least 3, and a vertex, respectively. A {\it dashed edge} and its two end vertices jointly represent an uncolored path of any length, and this path may be trivial. For example, in Figure~\ref{fig:9}, vertex $u$ is incident with the edge with color $a$ if the path is trivial. For any vertex, denoted by $q$, in the figures, uncolored $q$-branches may exist, even though they may not be shown in the figures. 

To facilitate the proof of the above theorem, we first introduce the concept of independent subtree. Throughout the game on $T$, colored edges can be interpreted as {\it boundaries} to decompose $T$ into {\it subtrees}, and each boundary belongs to exactly two subtrees. When an edge is being colored, the subtree containing this edge will be split into two subtrees. In particular, when a leaf-edge is colored, the subtree which contains it will be split into a $K_2$ and the subtree itself. Since it is clear that subsequent coloring of any subtree will not affect another, we may consider each subtree independently.

\begin{figure}
\centering
\subcaptionbox{An unpermitted 3-LCT \label{fig:312}} [.45\linewidth] 
 {\includegraphics[scale=.5]{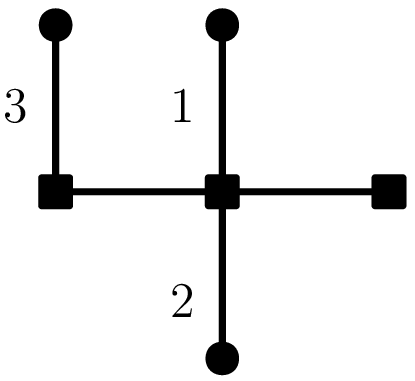}}
\subcaptionbox{An unpermitted 3-LCT \label{fig:212}} [.45\linewidth]
 {\includegraphics[scale=.5]{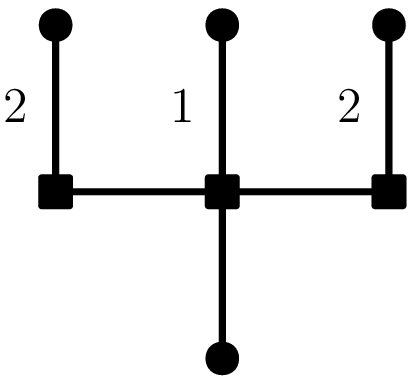}}
\subcaptionbox{An unpermitted 4-LCT \label{fig:1234}} [.45\linewidth]
 {\includegraphics[scale=.5]{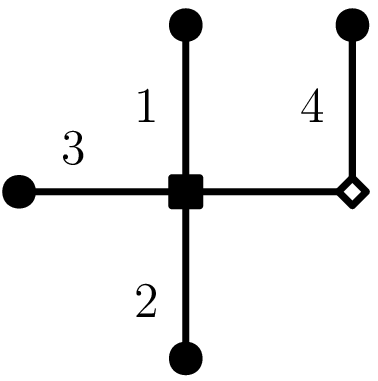}} 
\subcaptionbox{An unpermitted 4-LCT \label{fig:3123}} [.45\linewidth]
 {\includegraphics[scale=.5]{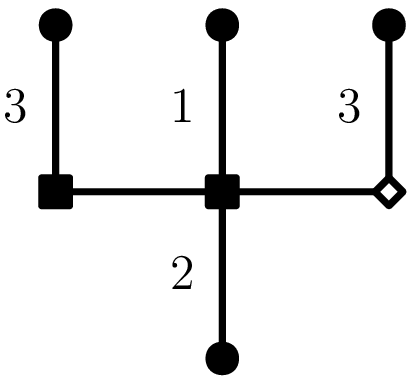}} 
\subcaptionbox{An unpermitted 4-LCT \label{fig:3124}} [.45\linewidth]
 {\includegraphics[scale=.5]{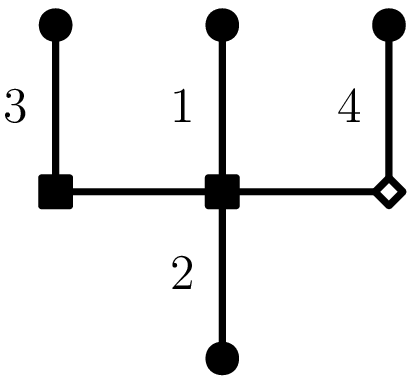}} 
\subcaptionbox{An unpermitted 4-LCT \label{fig:312a}} [.45\linewidth]
 {\includegraphics[scale=.5]{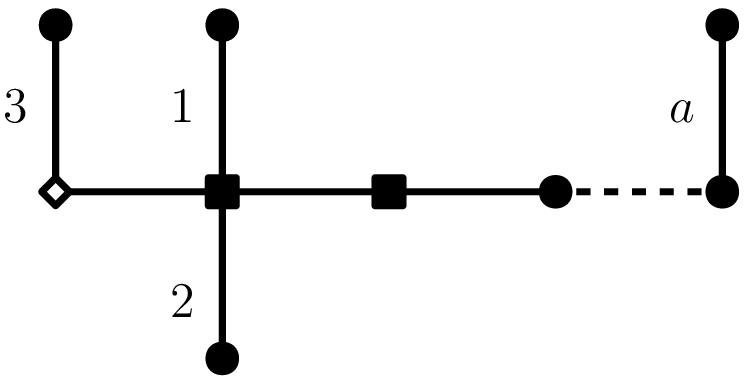}} 
\subcaptionbox{An unpermitted 4-LCT \label{fig:212a}} [.45\linewidth]
 {\includegraphics[scale=.5]{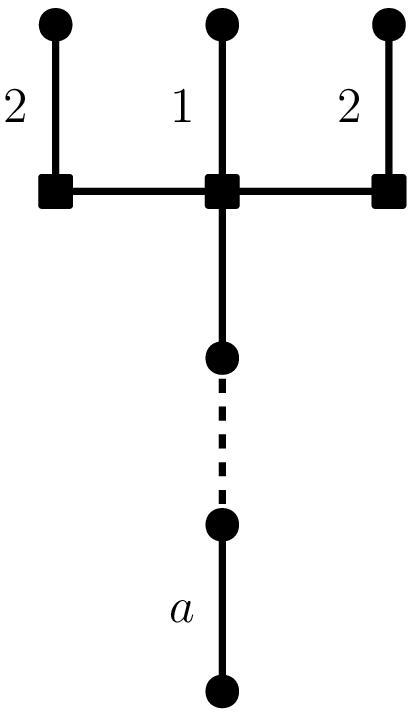}} 
	\caption{The unpermitted 3-LCTs and 4-LCTs with one star-node.}
\end{figure}

In all the figures and lemmas in this paper, our depictions and illustrations are based on colors of only some particular subset of the color set $\{1,2,3,4,5\}$, e.g., colors 1 and 2 in Figure~\ref{fig:212}, or colors 1, 2, 3 in Figure~\ref{fig:5}. Nevertheless, by the symmetry of the colors in the color set $\{1,2,3,4,5\}$, the depictions and illustrations based on a particular subset have no loss of generality so that they can be applied to cases with any other subset of the color set. As a result, the lemmas are not only true for some particular subset but also true for any other subset.

We call the following Types 0 to 40 subtrees {\it permitted}. Subtrees which are not permitted are called {\it unpermitted}.\\

\underline{ \textbf{Permitted types}}

\begin{figure}
\centering
\subcaptionbox{Type 11. $b\neq 1$.\label{fig:11}} [.45\linewidth]
 {\includegraphics[scale=.5]{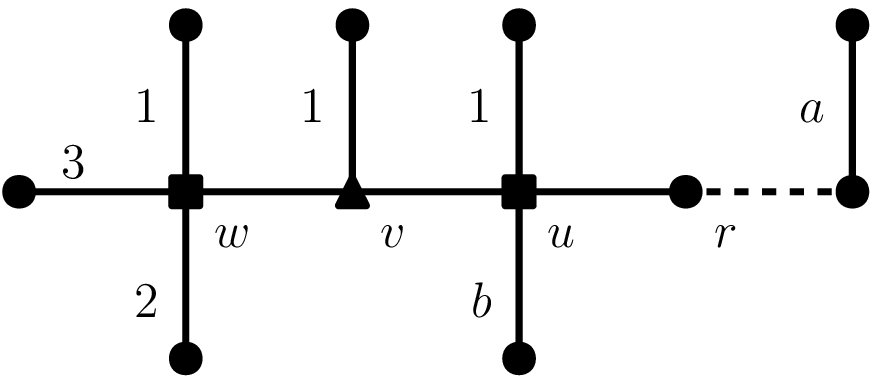}}  \hfill%
 \subcaptionbox{Type 12\label{fig:12}} [.45\linewidth]
 {\includegraphics[scale=.5]{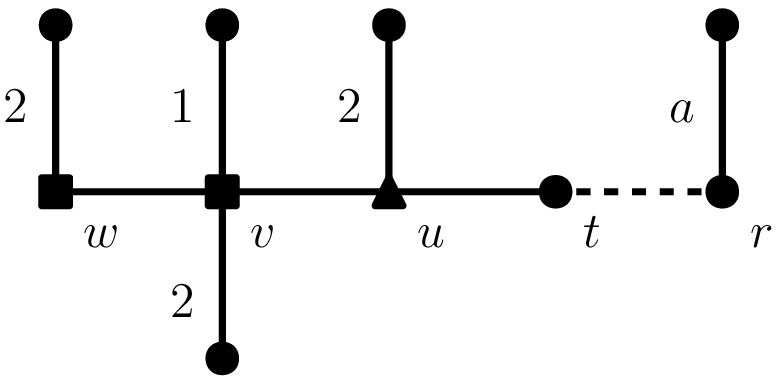}} 
\subcaptionbox{Type 13\label{fig:13}} [.45\linewidth]
 {\includegraphics[scale=.5]{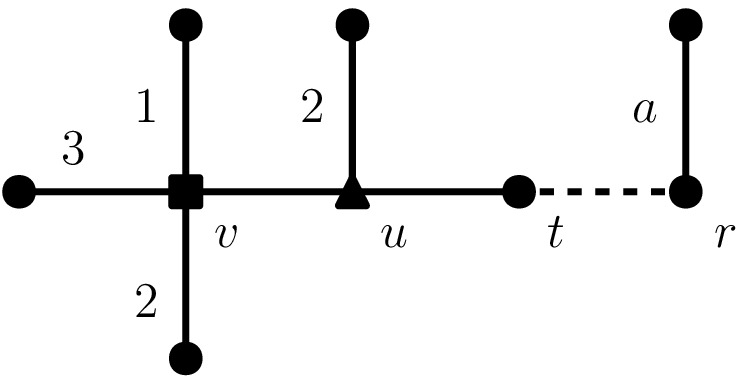}} \hfill%
\subcaptionbox{Type 14\label{fig:14}} [.45\linewidth]
 {\includegraphics[scale=.5]{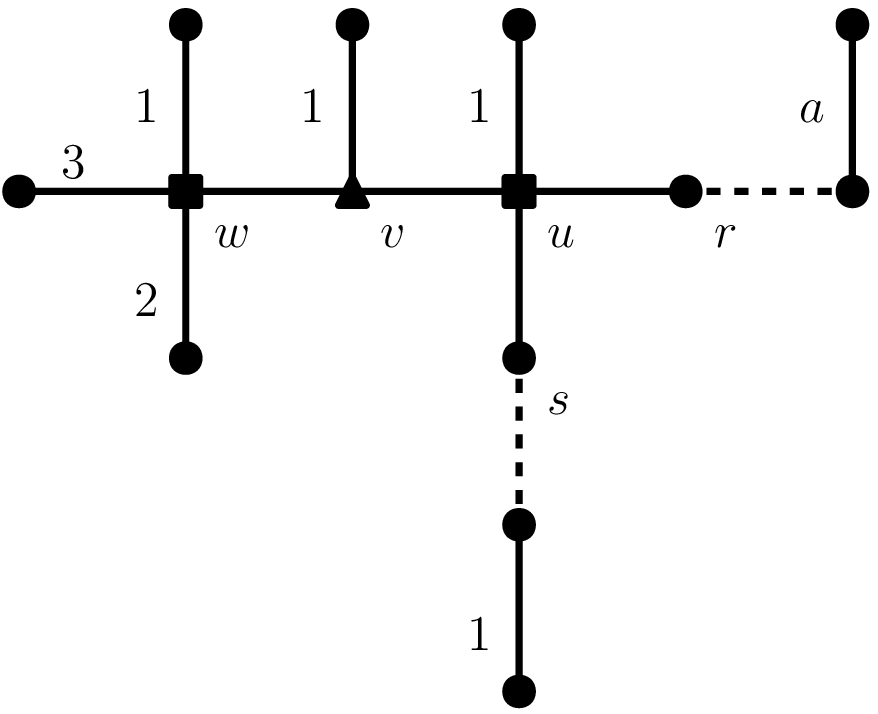}}  
\subcaptionbox{Type 15. $b\neq 1$.\label{fig:15}} [.45\linewidth]
 {\includegraphics[scale=.5]{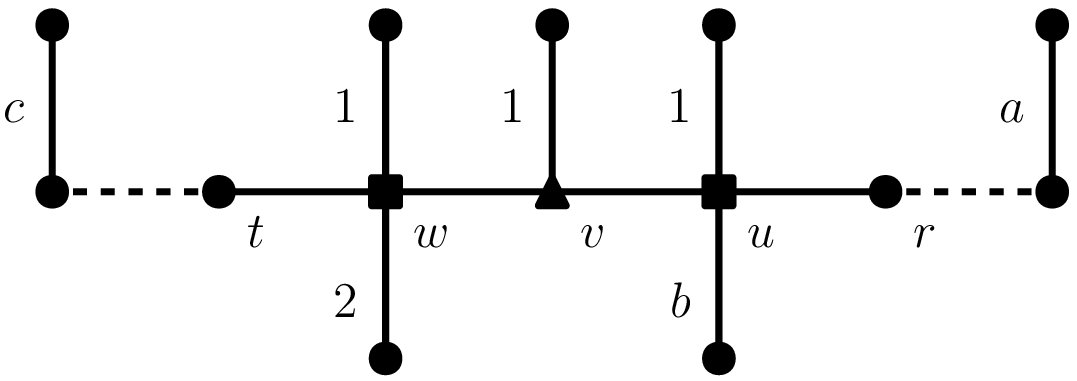}}  \hfill%
\subcaptionbox{Type 16\label{fig:16}} [.45\linewidth]
 {\includegraphics[scale=.5]{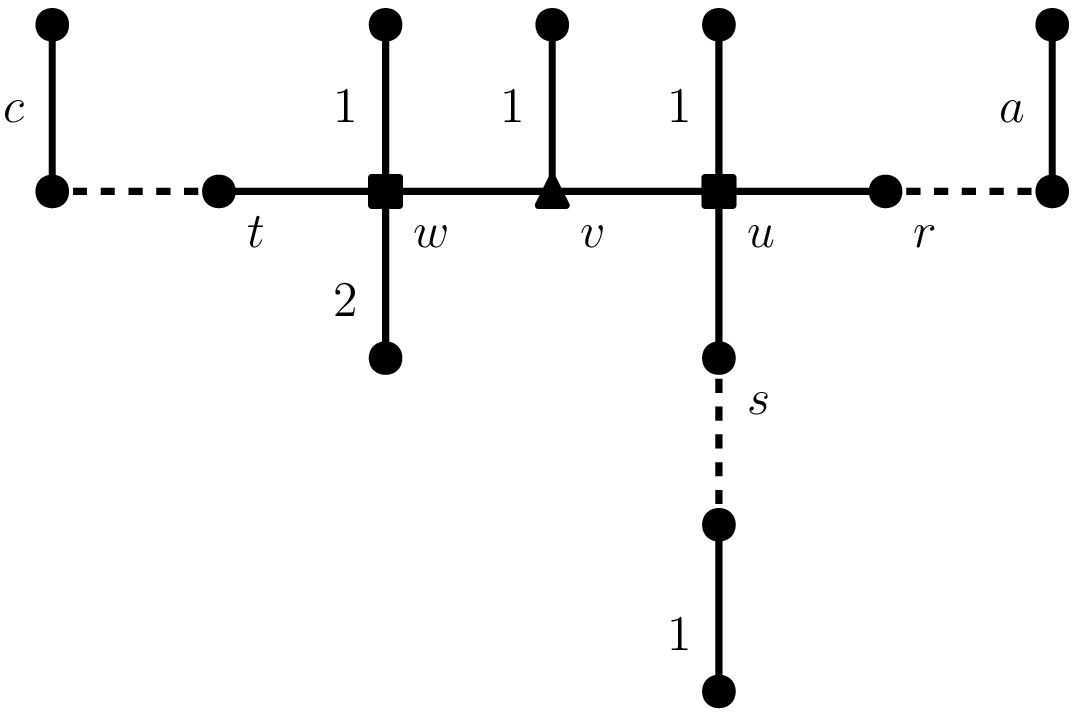}}  
 \subcaptionbox{Type 17. Vertex $y$ exists if $d(r)=4$. \label{fig:17}} [.45\linewidth]
 {\includegraphics[scale=.5]{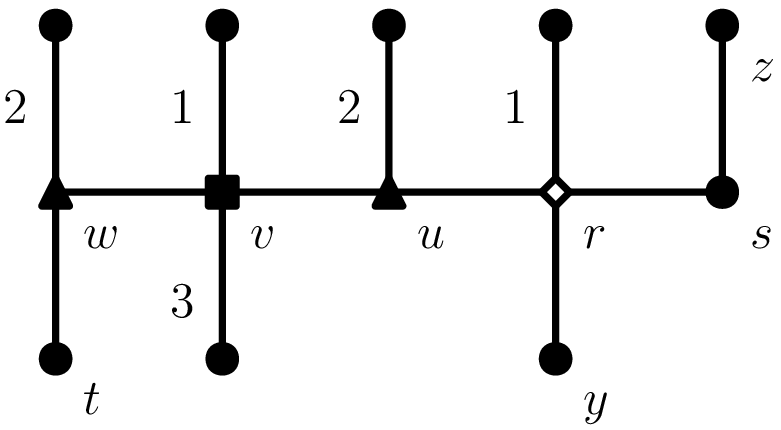}} \hfill%
\subcaptionbox{Type 18\label{fig:18}} [.45\linewidth]
 {\includegraphics[scale=.5]{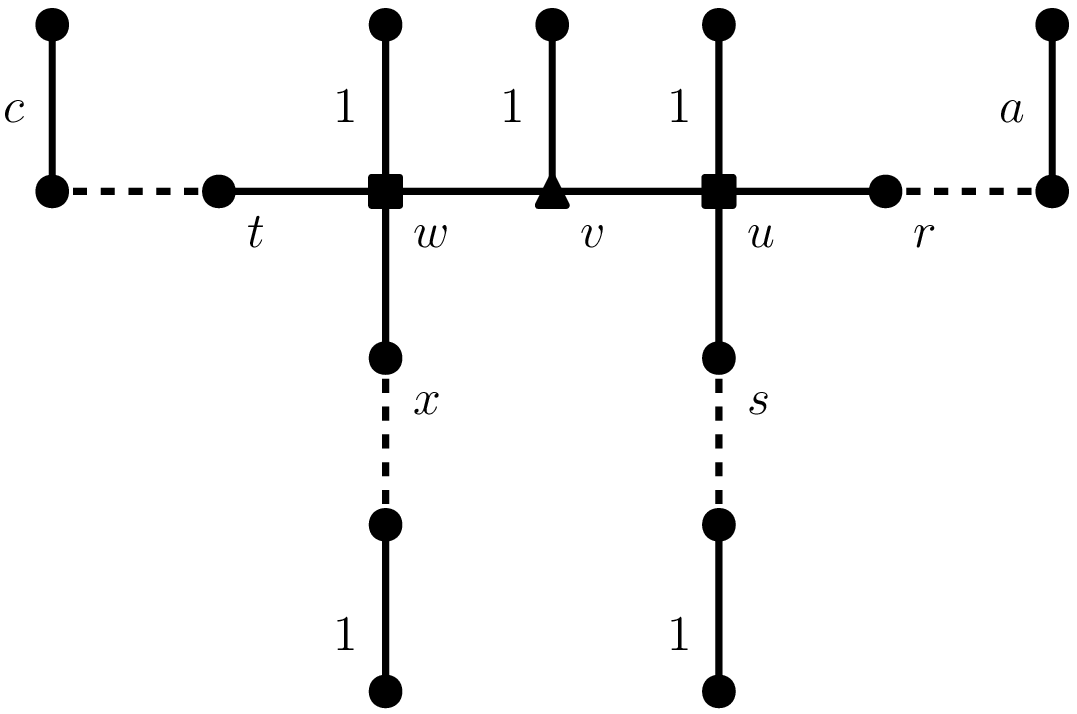}} 
\subcaptionbox{Type 19. It is under the trivial constraint that the two colored edges incident with $u$ do not have the same color. \label{fig:19}} [.45\linewidth]
 {\includegraphics[scale=.5]{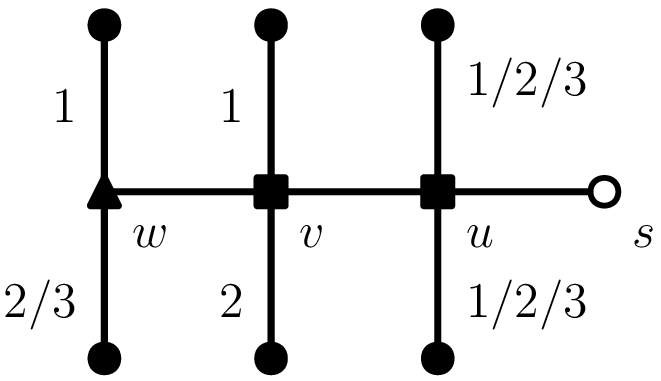}} \hfill%
\subcaptionbox{Type 20 \label{fig:20}} [.45\linewidth]
 {\includegraphics[scale=.5]{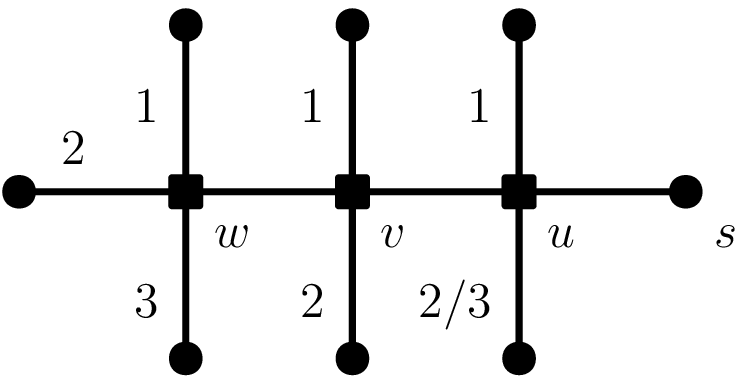}} 
	\caption{Types 11 - 20.\label{fig:11-20}}
\end{figure}

\begin{figure}
\centering

\subcaptionbox{Type 21. Since $S_3^4$ and $P_4^4$ are forbidden, $d(r),d(t)\leq 3$. Vertex $x$ exists if $d(r)=3$. \label{fig:21}} [.45\linewidth]
 {\includegraphics[scale=.5]{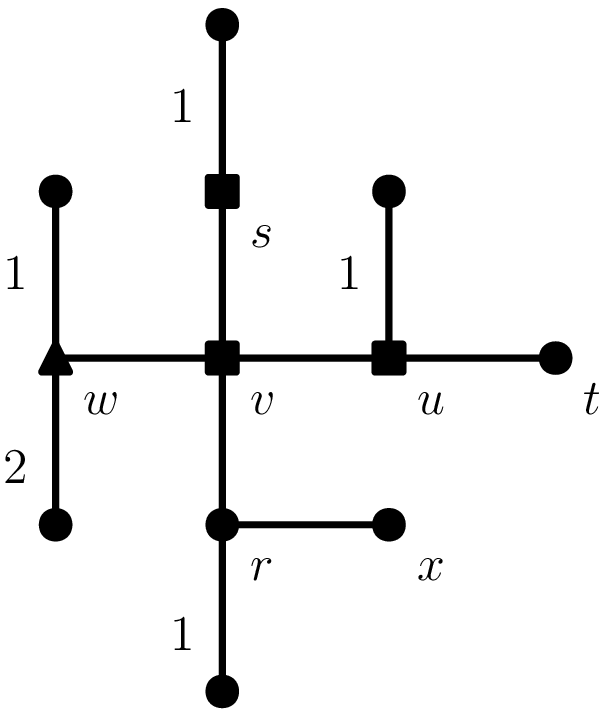}} \hfill%
\subcaptionbox{Type 22. Vertex $y$ exists if $d(r)=4$. \label{fig:22}} [.45\linewidth]
 {\includegraphics[scale=.5]{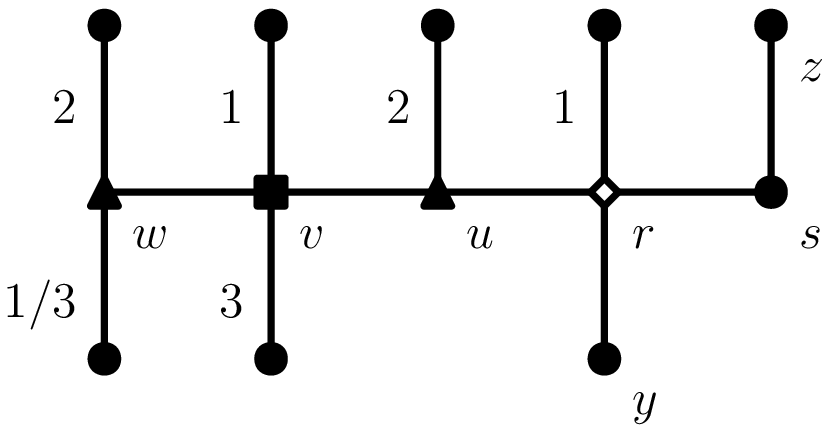}} 
 \subcaptionbox{Type 23  \label{fig:23}} [.45\linewidth]
  {\includegraphics[scale=.5]{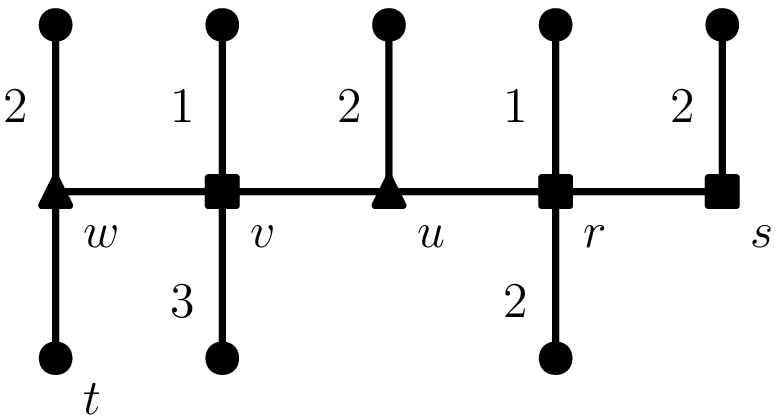}} \hfill%
 \subcaptionbox{Type 24  \label{fig:24}} [.45\linewidth]
 {\includegraphics[scale=.5]{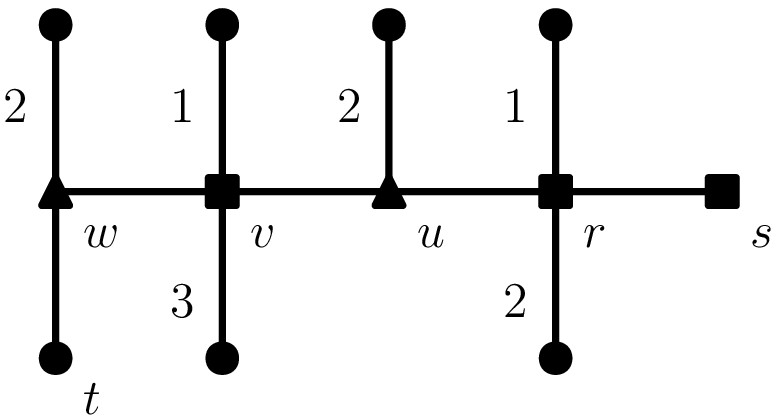}} 
 \subcaptionbox{Type 25  \label{fig:25}} [.45\linewidth]
 {\includegraphics[scale=.5]{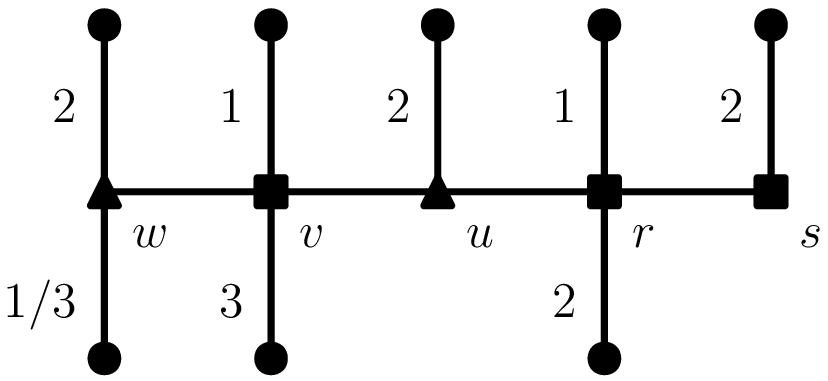}} \hfill%
 \subcaptionbox{Type 26 \label{fig:26}} [.45\linewidth]
 {\includegraphics[scale=.5]{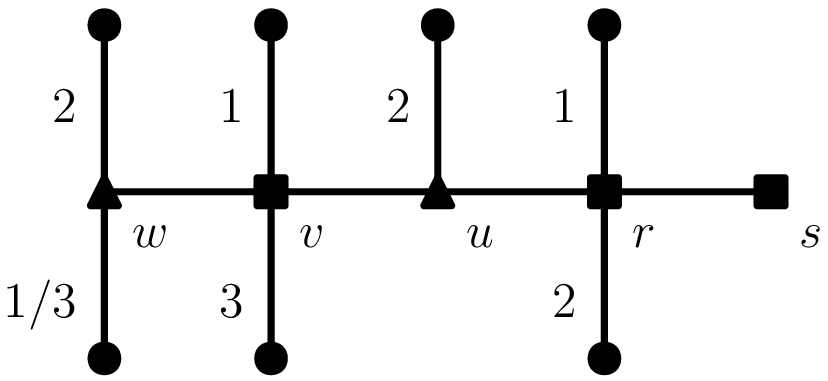}} 
  \subcaptionbox{Type 27 \label{fig:27}} [.45\linewidth]
 {\includegraphics[scale=.5]{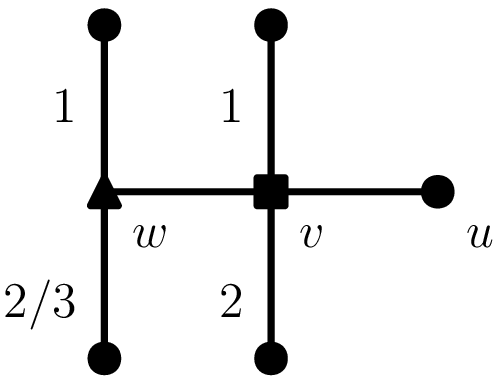}} \hfill%
   \subcaptionbox{Type 28 \label{fig:28}} [.45\linewidth]
 {\includegraphics[scale=.5]{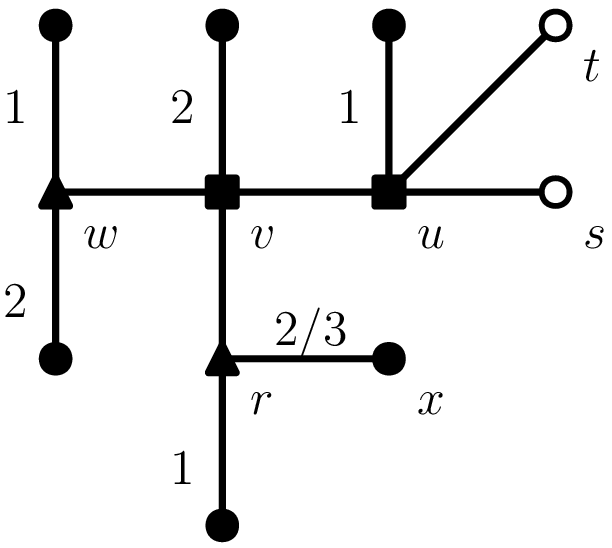}} 
    \subcaptionbox{Type 29 \label{fig:29}} [.45\linewidth]
 {\includegraphics[scale=.5]{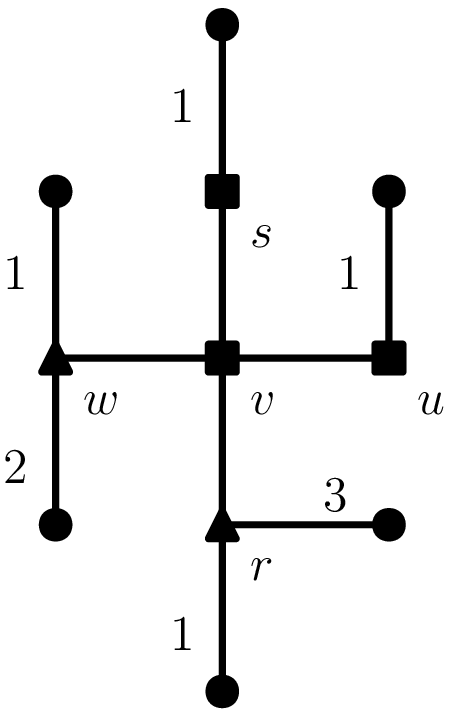}} \hfill%
     \subcaptionbox{Type 30. Since $S_3^4$ and $P_4^4$ are forbidden, $d(r),d(s),d(t)\leq 3$. Vertex $x$ exists if $d(r)=3$. \label{fig:30}} [.45\linewidth]
 {\includegraphics[scale=.5]{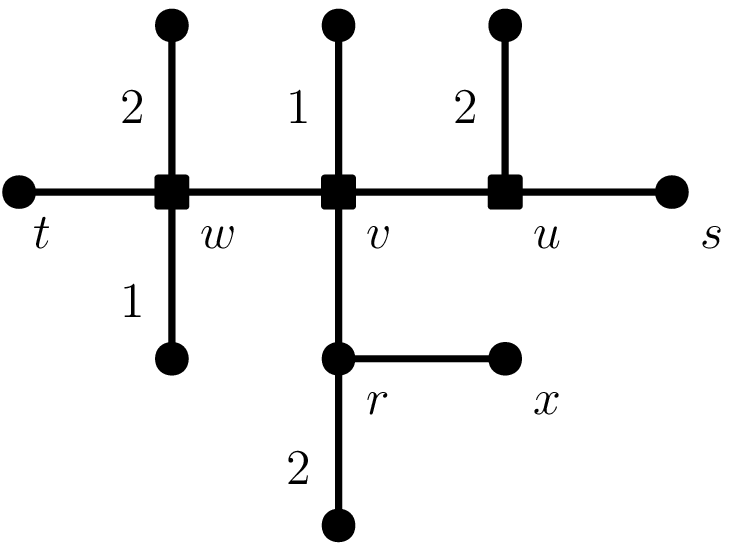}} 
	\caption{Types 21 - 30.\label{fig:21-30}}
\end{figure}

\begin{figure}
\centering
      \subcaptionbox{Type 31 \label{fig:31}} [.45\linewidth]
 {\includegraphics[scale=.5]{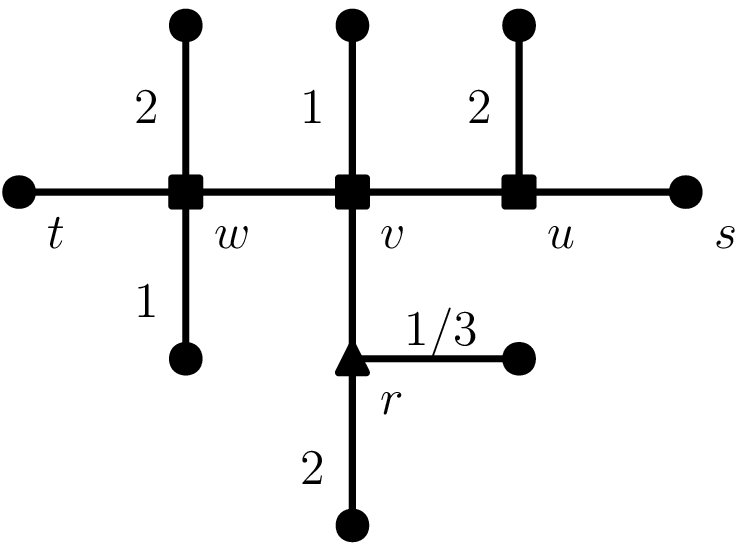}} \hfill%
\subcaptionbox{Type 32. $a=1$ or 3. \label{fig:32}} [.45\linewidth]
 {\includegraphics[scale=.5]{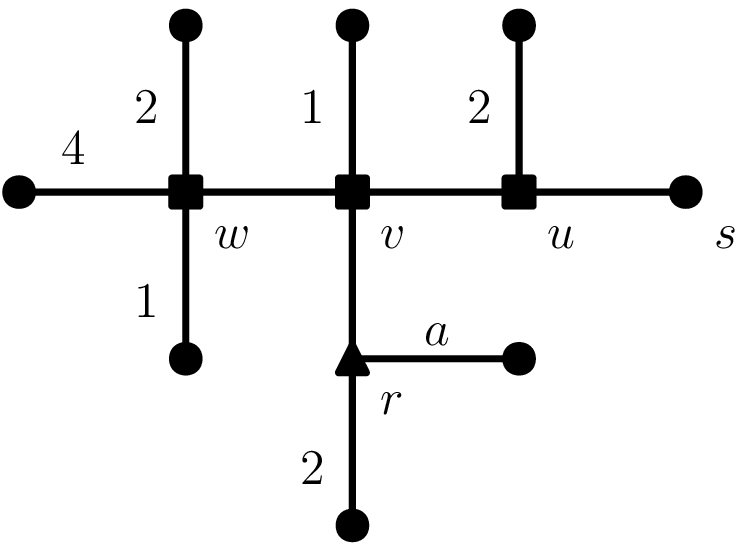}} 
\subcaptionbox{Type 33 \label{fig:33}} [.45\linewidth]
 {\includegraphics[scale=.5]{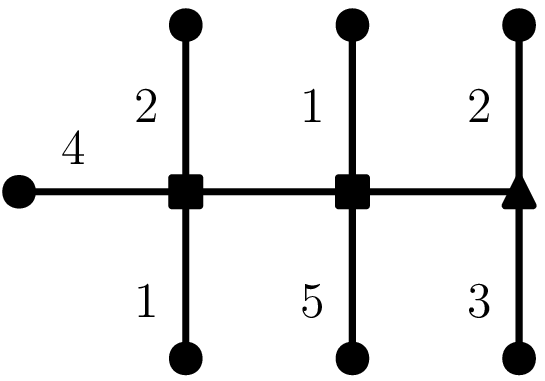}} \hfill%
\subcaptionbox{Type 34. Since $S_3^4$ and $P_4^4$ are forbidden, $d(r),d(s),d(t)\leq 3$. Vertex $x$ exists if $d(r)=3$. \label{fig:34}} [.45\linewidth]
 {\includegraphics[scale=.5]{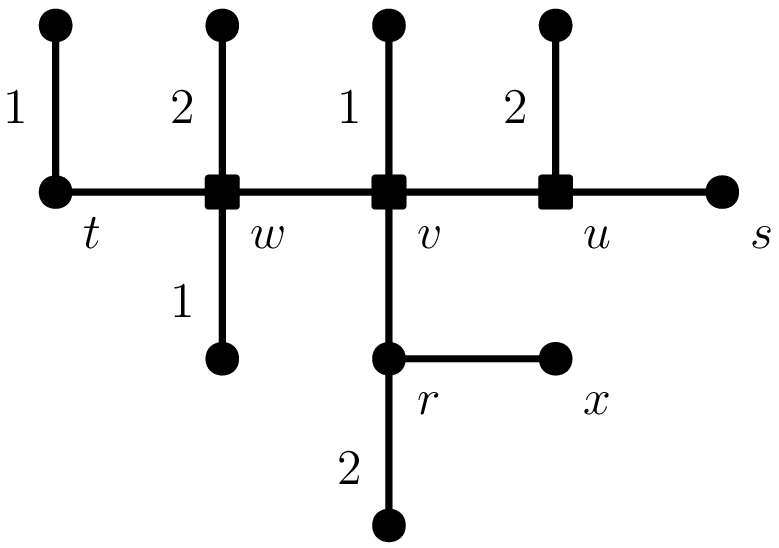}} 
\subcaptionbox{Type 35 \label{fig:35}} [.45\linewidth]
 {\includegraphics[scale=.5]{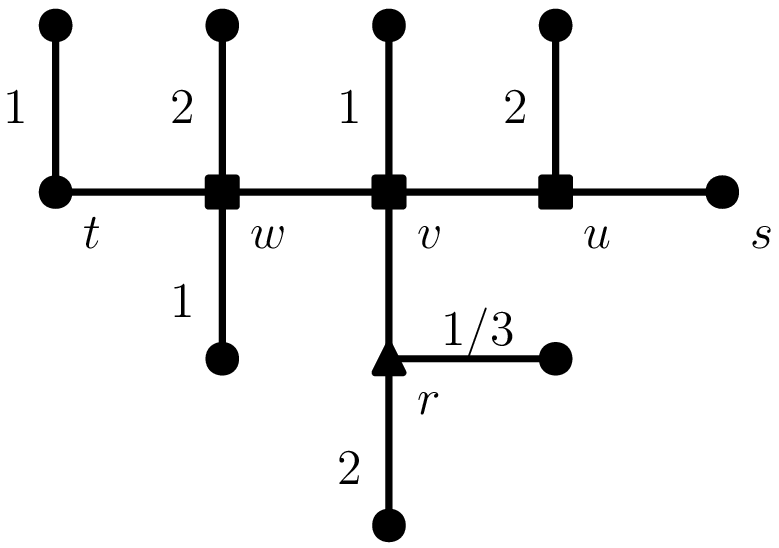}} \hfill%
\subcaptionbox{Type 36. Vertex $s$ (resp. $t$) exists if $d(r)=3$ (resp. $d(u)=3$). \label{fig:36}} [.45\linewidth]
 {\includegraphics[scale=.5]{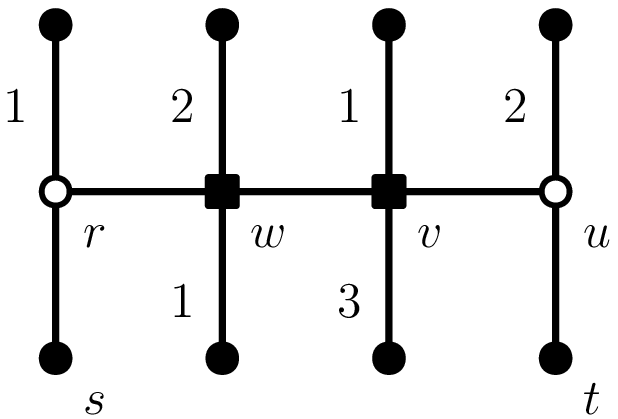}} 
 \subcaptionbox{Type 37 \label{fig:37}} [.45\linewidth]
 {\includegraphics[scale=.5]{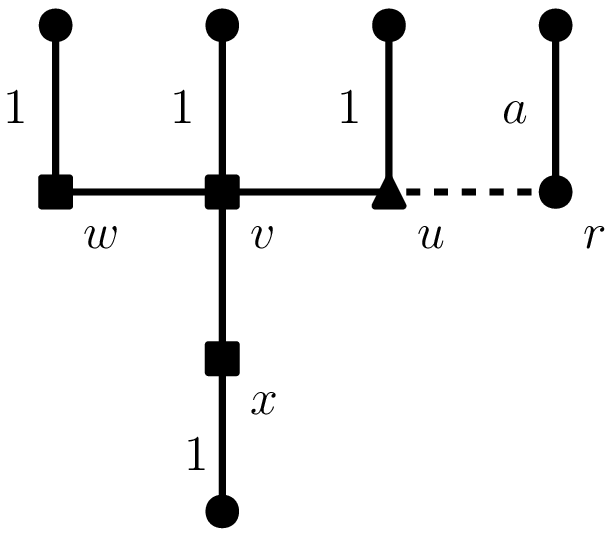}} \hfill%
 \subcaptionbox{Type 38. It is under the constriant that if $d(u)=4$, then vertex $s$ exists and $d(t),d(s)\leq 3$. \label{fig:38}} [.45\linewidth]
 {\includegraphics[scale=.5]{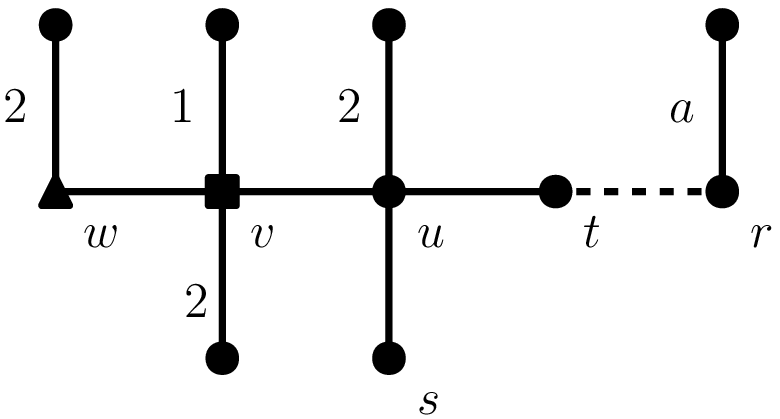}} 
 \subcaptionbox{Type 39. Vertex $x$ exists if $d(r)=3$. \label{fig:39}} [.45\linewidth]
 {\includegraphics[scale=.5]{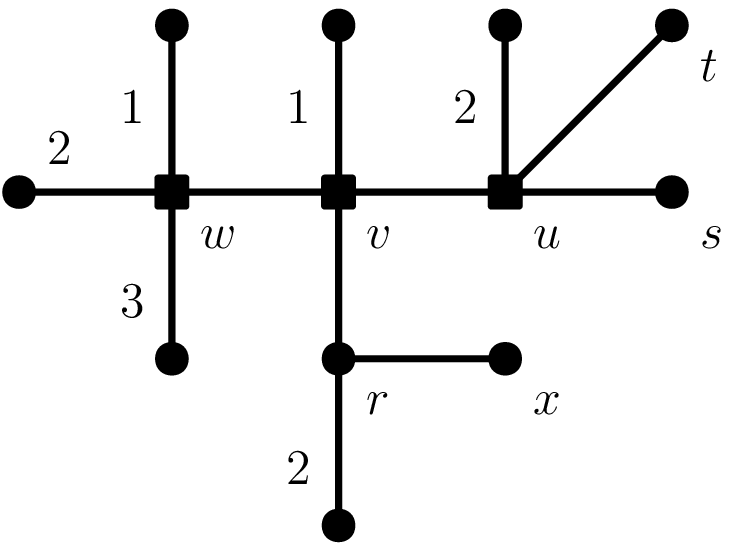}} \hfill%
 \subcaptionbox{Type 40\label{fig:40}} [.45\linewidth]
 {\includegraphics[scale=.5]{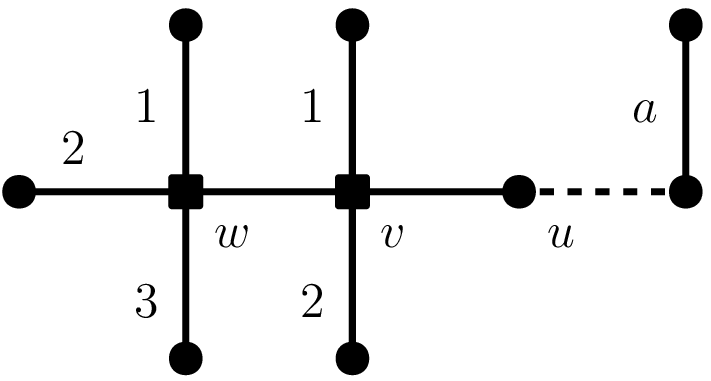}} 
	\caption{Types 31 - 40.\label{fig:31-40}}
\end{figure}

\begin{itemize}
\item Type 0: Completely colored subtrees. 
\item Type 1: Subtrees containing no star-nodes. In particular, completely colored stars with at most two leaves are of Type 1.
\item Type 2: Subtrees containing exactly one star-node, except those in Figures~\ref{fig:312} to~\ref{fig:212a}. In particular, Completely colored stars with three or four leaves are of Type 2.
 
We remark that:
\begin{enumerate}[i.] 
\item The subtrees in Figures~\ref{fig:312} and~\ref{fig:212} are 3-LCT, and those in Figures~\ref{fig:1234},~\ref{fig:3123},~\ref{fig:3124},~\ref{fig:312a} and~\ref{fig:212a} are 4-LCT. 
\item Each colored edge of Figures~\ref{fig:312} and~\ref{fig:212} is incident with some 4-vertex. 
\item The subtree in Figures~\ref{fig:1234} has three colored star-edges; Figures~\ref{fig:3123},~\ref{fig:3124} and~\ref{fig:312a} have two colored star-edges; and Figure~\ref{fig:212a} has one colored star-edge.
\item In the subtrees in Figures~\ref{fig:312} -~\ref{fig:312a} but not in Figure~\ref{fig:212a}, any colored star-edge does not share the same color with any colored edge which is incident with a neighbor of the star-node.
\end{enumerate}  
The above properties will be employed in the rest of this paper.

\item Type 3: A subtree with the subgraph induced by all its uncolored edges being a path $P$ of length $m$, where $1\leq m\leq 3$. Moreover, for $m>1$, at least two colors are available for each uncolored edge; the only uncolored edge has at least one available color for $m=1$.
\item Type 4: A subtree with the subgraph induced by all its uncolored edges being the union of a path $P=v_0v_1...v_m$, where $3\leq m\leq 4$, and a tree $T_s$ (maybe trivial). Also, this subtree has to satisfy the two requirements that (1) $v_m$ is the only common vertex of $P$ and $T_s$, and (2) all vertices of $T_s$ are not incident with any colored edge. Moreover, the path $P$ has to satisfy the condition that $v_{m-1}v_m$ has at least four available colors and each of the other uncolored edges of $P$ has at least two available colors.  An example of this type with $m=4$ is shown in Figure~\ref{fig:4ex}.
\item Types 5-40: See Figures~\ref{fig:4-10},~\ref{fig:11-20},~\ref{fig:21-30}, and~\ref{fig:31-40}.
\end{itemize}

%


\bigskip

To show that the tree $T$ is 5-edge-game-colorable, we will propose a strategy of Alice for picking an edge and assigning a color on it in each of her turns in the proof of Lemma~\ref{lem:inductive}. Before that, we first show that Theorem~\ref{mainthm} can be derived from Lemma~\ref{lem:inductive}.  Meanwhile, we give another lemma which ensures that after any Bob's move on a permitted subtree, at most one generated subtree is unpermitted.

\begin{lemma}\label{lem:inductive}
If all subtrees are permitted after a move of Alice, then Alice may keep all subtrees permitted right after her next move regardless of the move of Bob in this turn.
\end{lemma}

\medskip

\noindent {\bf Proof of Theorem~\ref{mainthm}.} Under the game variant that Bob can skip, Alice loses the game on a finite tree if and only if she cannot find any available move in some of her turns. As a result, it suffices to show that Alice has a feasible move in each of her turns. 

We will show that in each of Alice's turns, she can always find a move such that right after this move, all subtrees are permitted. This can be proved by induction on the number of Alice's turns. In her first turn, any subtree that Alice receives is either a 0- or 1-LCT. Hence right after her first move, each subtree has less than three colored edges, which implies all subtrees are of Type 0 or 1. Therefore, if Lemma~\ref{lem:inductive} holds, by induction, she can ensure that all subtrees are permitted right after anyone of her moves. \eop  

\begin{lemma}\label{lem:1permitted}
Suppose Bob colors an edge of a permitted subtree $F$ and generates two subtrees. At most one of them is unpermitted.
\end{lemma}

\medskip

\noindent {\bf Proof of Lemma~\ref{lem:1permitted}.} First, any subtree containing at most two colored edges is permitted; therefore, at least six edges have been colored in the two subtrees in total if they are both unpermitted. Second, after Bob has colored an edge $e$ of $F$, the two newly generated subtrees should have two more colored edges than that were in $F$, because $e$ is double counted in the two subtrees. Consequently, the lemma holds if the number of colored edges of $F$ is at most 3. We remark that any Type 1 subtree and any 3-LCT of Type 2 has at most three colored edges. In the following, we will deal with the rest of the permitted types.

\begin{itemize}
\item Suppose $F$ is a 4-LCT of Type 2. A 4-LCT will be decomposed into two 3-LCT only if an edge on the path connecting two star-nodes is colored. As $F$ has exactly one star-node, at least one of the newly formed subtrees has less than three colored edges, which is of Type 1.

\item Suppose $F$ is of Type 3 or 33. Bob can only color any edge of the uncolored path of length at most three.  At least one of the subtrees formed is of Type 0 or 3.

\item Suppose $F$ is of Type 4. When Bob colors any edge of the uncolored path $P$, one of the two newly formed subtrees must be of Type 0, 1 or 3; when he colors any other edge, one of the two newly formed subtrees must contain only one colored edge.

\item  Suppose $F$ is of Type 15, 16 or 18. 

If Bob colors any edge incident with a colored edge, except $wv$ and $vu$, a Type 1 or 2 subtree will be generated. He will generate a Type 1 subtree when he colors any edge not incident with a colored edge. 

If he colors $wv$, the generated subtree which does not contain vertex $t$ is of Type 9 (when $F$ is of Type 15) or Type 10 (when $F$ is of Type 16/18). 

If he colors $vu$, the generated subtree which does not contain vertex $r$ is of Type 9 (when $F$ is of Type 15/16) or Type 10 (when $F$ is of Type 18).

\item Suppose $F$ is of Type 23 or 24. If Bob colors any edge incident with a colored edge, except $vu$, a Type 1 or Type 2 subtree will be generated. He will generate a Type 1 subtree when he colors an edge not incident with any colored edge. 

If he colors $vu$, the generated subtree which does not contain $w$ is of Type 9.

\item Suppose $F$ is of Type 5, 7, 8, 9, 10, 11, 13, 14, 19, 20, 21, 22, 27, 28, 29, 32, 39 or 40. If Bob colors any edge incident with a colored edge, a Type 0/1/3 subtree will be generated. He will generate a Type 1 subtree when he colors an edge not incident with any colored edge. 

\item Suppose $F$ is of Type 6, 12, 17, 25, 26, 30, 31, 34, 35, 36, 37 or 38. If Bob colors any edge incident with a colored edge, a Type 0/1/2/3 subtree will be generated. He will generate a Type 1 subtree when he colors an edge not incident with any colored edge. \eop
\end{itemize}

In the next two sections, we will prove Lemma~\ref{lem:inductive} using Lemma~\ref{lem:1permitted}.

\section{Proof of Lemma~\ref{lem:inductive}: Alice's strategy for Types 1 and 2 subtrees}
\medskip
In Lemma~\ref{lem:1permitted}, we have shown that Bob generates at most one unpermitted subtree in his act on a permitted subtree. Therefore, in each turn of Alice, in order to ensure that all subtrees are permitted right after her move, her task is to turn $F$, the unpermitted subtree if any, or, otherwise, a non-completely colored permitted subtree, to two permitted subtrees. In the rest of this section, we shall propose a strategy of Alice for handling (i) subtrees of Types 1 and 2; and (ii) the unpermitted subtrees which would be made by any possible move of Bob on Types 1 and 2 subtrees. In the next section, we will introduce, one by one, strategies of Alice for handling (i) Types 3 to 40 subtrees; and (ii) the unpermitted subtrees which would be made by any possible move of Bob on Types 3 to 40 subtrees, respectively. Note that 1) all non-completely colored subtrees are permitted just before any move of Bob; and 2) a non-completely colored subtree may be a member of two different permitted types simultaneously, for example, Types 2 and 3. For this case, Alice may employ any one of the two strategies for the two distinct types of subtrees.

\medskip
First, note that for case (i) (that is, handling Types 1 and 2 subtrees), all Types 1 and 2 subtrees have zero and one star-node, respectively. Second, we consider case (ii) (that is, handling the unpermitted subtrees which would be generated by any move of Bob on Types 1 and 2 subtrees). For those generated from a Type 1 subtree, at most one of the two generated subtrees contains a star-node. Any unpermitted subtree obtained from Type 2, which contains at most four colored edges and precisely one star-node, has three, four or five colored edges. Moreover, any 5-LCT generated by coloring an edge of a permitted 4-LCT (which means this 4-LCT has one 4-SN) has precisely one 4-SN and one 3-SN. To conclude, $F$, the subtree on which Alice is going to act, has zero, one or two star-nodes. In the rest of this section, we will propose a unified strategy of Alice for handling cases (i) and (ii) together with respect to the three situations of the number of star-nodes in $F$. Unless specified otherwise, Alice may make use of any available color that can generate subtrees of the desired types.

\noindent (I) A non-completely colored $F$ with no star-nodes:
\bit
\item $F$ is a 0-LCT or a 1-LCT:

Alice may color any edge to generate two Type 1 subtrees.

\item $F$ is a 2-LCT:

If the two colored edges are adjacent, Alice may color an edge adjacent to both these two colored edges to generate one Type 1 and one Type 2 subtree. 

If the two colored edges are not adjacent, Alice may color an edge on the path containing the two roots of the colored edges to generate two Type 1 subtrees.
\eit

\noindent (II) A non-completely colored 3-LCT $F$ with exactly one star-node:
\begin{itemize}
\item If $F$ has three colored star-edges, Alice may color the remaining star-edge to make a Type 0 and a Type 1 subtree.

\item If $F$ has exactly two colored star-edges, Alice may color the star-edge on the leaf-path to make one Type 1 and a Type 0 (if $d(v)=3$) or Type 2 (if $d(v)=4$) subtree.

\item Suppose $F$ has exactly one colored star-edge. If $d(v)=3$, Alice may color any star-edge to make one Type 1 and one 3-LCT. This 3-LCT is of Type 2 since an unpermitted 3-LCT must have a SN of degree 4. 

When $d(v)=4$, we denote the two vertices adjacent to $v$ and on the two leaf-paths by $u$ and $w$. We also denote the vertex adjacent to $v$ but not on any leaf-path by $t$. When at least one of vertices $u$ and $w$ is a 4-vertex, without loss of generality, we assume $d(u)=4$. Then, Alice may color $vu$ to generate a Type 1 subtree and a 3-LCT.  Given that $S_3^4$ is forbidden, the star-node $v$ of this 3-LCT is adjacent to at most one 4-vertex in this 3-LCT, which means this 3-LCT is not the one in Figure~\ref{fig:312} nor~\ref{fig:212}. This implies the 3-LCT is of Type 2. When both $u$ and $w$ are not of degree 4, Alice may color $vt$ to generate a Type 1 subtree and a 4-LCT with 2 colored star-edges. This 4-LCT is of Type 2 since its SN $v$ is not adjacent to any 4-vertex in this 4-LCT.

\item Suppose $F$ has no colored star-edge, Alice may use a strategy similar to that for the case that $F$ having exactly one colored star-edge: 

If $d(v)=3$, Alice may color any star-edge to make one Type 1 and one 3-LCT of Type 2.

If $d(v)=4$, we denote the three vertices adjacent to $v$ and on the three leaf-paths by $u$, $w$ and $r$. $t$ denotes the vertex adjacent to $v$ but not on any leaf-path. When at least one of $u$, $w$ and $r$ is a 4-vertex, without loss of generality, we assume $d(u)=4$. Then, Alice may color $vu$ to generate a Type 1 subtree and a 3-LCT of Type 2. When all $u, w$ and $r$ are not of degree 4, Alice may color $vt$ to generate a Type 1 subtree and a 4-LCT with 1 colored star-edge and of Type 2. 
\end{itemize}

\noindent (III) A non-completely colored 4-LCT $F$ with exactly one star-node:

\begin{itemize}
\item If $F$ has exactly three colored star-edges, Alice may color the remaining star-edge to make one Type 0 and one Type 1 subtree.

\begin{figure}
\centering
\subcaptionbox{$a$ and $b$ can be arbitrary colors. \label{fig:4lct2cse}} [.45\linewidth]
 {\includegraphics[scale=.5]{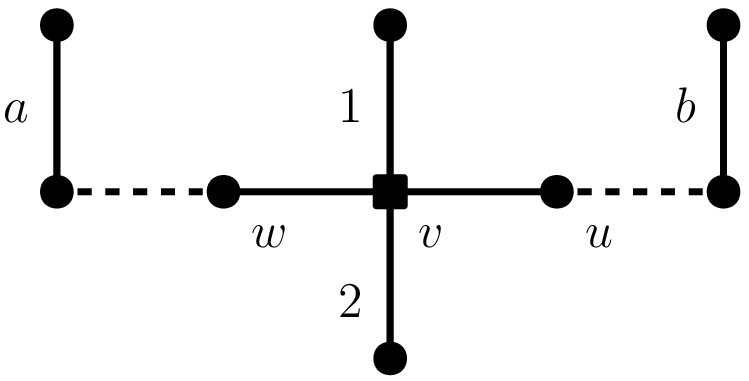}} 
 \subcaptionbox{Note that $d(u)\geq 3$. \label{fig:3123d33}} [.45\linewidth]
 {\includegraphics[scale=.5]{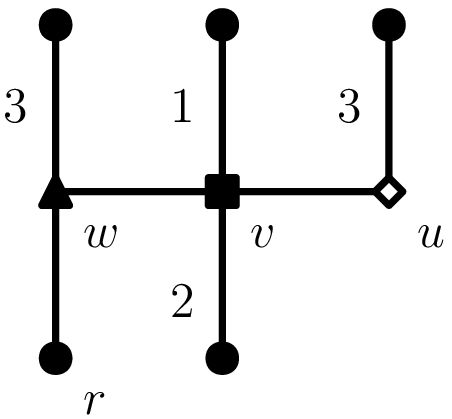}} 
	\caption{Some subtrees $F$ having exactly two colored star-edges}
\end{figure}

\item Suppose $F$ has exactly two colored star-edges, then $F$ can be represented by the subtree in Figure~\ref{fig:4lct2cse}. Note that this $F$ is symmetric. Also, only colors 1 and 2 are used in the figure so colors 3, 4 and 5 are identical by symmetry.

If $a=1$ or 2 (resp. $b=1$ or 2), Alice may color $vu$ (resp. $vw$) with 3/4/5 to generate one Type 1 subtree and one 4-LCT $T_v$ with three colored star-edges. The four colored edges of $T_v$ have three different colors only, which means it is not the one in Figure~\ref{fig:1234}; so $T_v$ is of Type 2.

Suppose both $a$ and $b$ are not in $\{1,2\}$. Alice may color $vu$ by $a$ or $vw$ by $b$, if feasible, to generate one Type 1 subtree and one 4-LCT of Type 2 with three colored star-edges. Such two acts of Alice are both impossible if only if $a=b=3$ and the edges with colors $a$ and $b$ are incident with $w$ and $u$, respectively. We then discuss alternative strategies of Alice based on the degrees of $w$ and $u$.

Suppose $d(u)=d(w)=4$, then $F$ is the unpermitted 4-LCT in Figure~\ref{fig:3123}, which means this $F$ was just generated from a permitted subtree by Bob so some edge of this $F$ was just colored by Bob. If Bob just colored the edge with color $b=3$, $F$ was generated either from the unpermitted 3-LCT in Figure~\ref{fig:312} or the 4-LCT in Figure~\ref{fig:312a}, which leads to a contradiction. If Bob just colored the edge with color 1 or 2, $F$ was generated either from the unpermitted 3-LCT in Figure~\ref{fig:212} or the 4-LCT in Figure~\ref{fig:212a}, which are both impossible. To conclude, at most one of $u$ and $w$ is a 4-vertex. We will discuss it with two cases: (1) $d(w),d(u)\geq 3$; or (2) at least one of $w$ and $u$ is of degree 2. 

(1) Without loss of generality, assume $d(w)=3$ and $d(u)\geq 3$; hence, $F$ can be represented by the one in Figure~\ref{fig:3123d33}. Alice may color $rw$ with 1 to make one Type 1 and one Type 9 subtree. 

(2) Without loss of generality, assume $d(w)=2$. Alice may color $vu$ with color 4 to generate a Type 3 and a Type 1 subtree.

\item Suppose $F$ has exactly one colored star-edge. If $v$ is adjacent to exactly one 4-vertex, denoted by $u$, Alice may color $vu$ to generate one Type 1 subtree and one 4-LCT with exactly two colored star-edges. This 4-LCT is of Type 2 since its SN $v$ is not adjacent to any 4-vertex. Similarly, if $v$ is not adjacent to any 4-vertex, Alice may color any star-edge to generate one Type 1 subtree and one 4-LCT of Type 2 with exactly two colored star-edges.

\begin{figure}
\centering
\subcaptionbox{Since $S_3^4$ is forbidden, $d(r)\leq 3$. \label{fig:4lct1cse444}} [.45\linewidth]
 {\includegraphics[scale=.5]{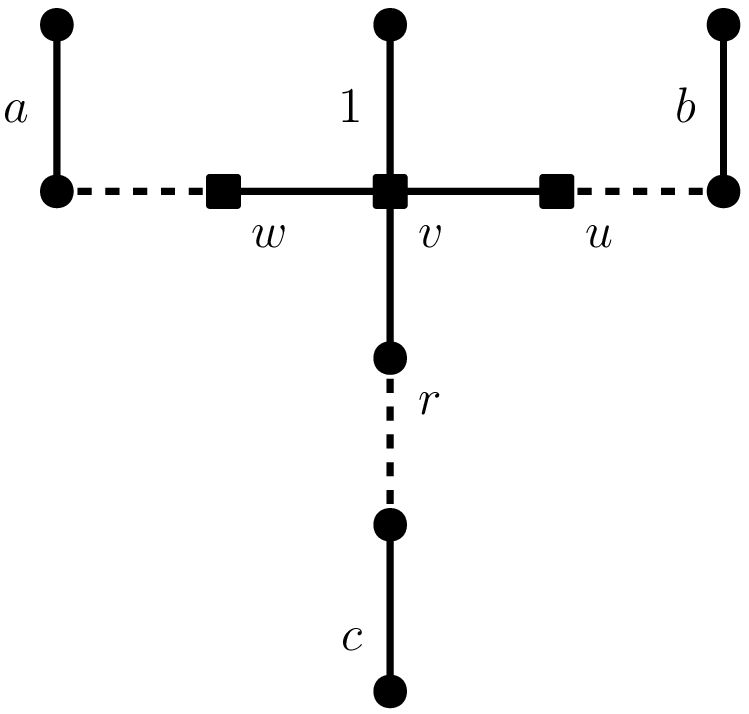}} \hfill%
\subcaptionbox{Coloring $vu$ with 2 would lead to an unpermitted subtree. \label{fig:4lct1cse3124}} [.45\linewidth]
 {\includegraphics[scale=.5]{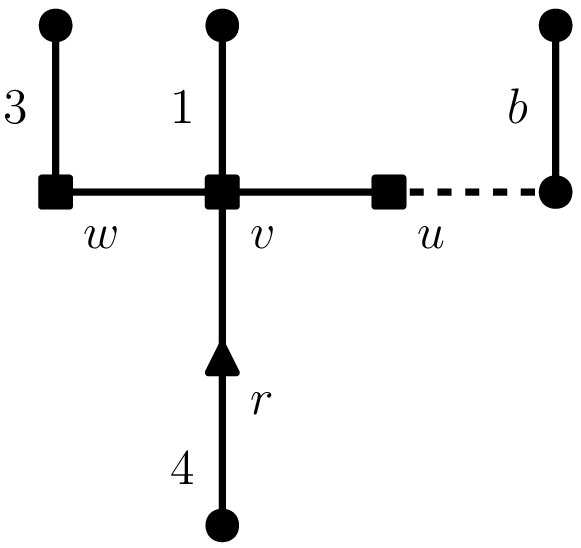}}
 
\subcaptionbox{Coloring $vu$ with 2 would lead to an unpermitted subtree. \label{fig:4lct1cse3123}} [.45\linewidth]
 {\includegraphics[scale=.5]{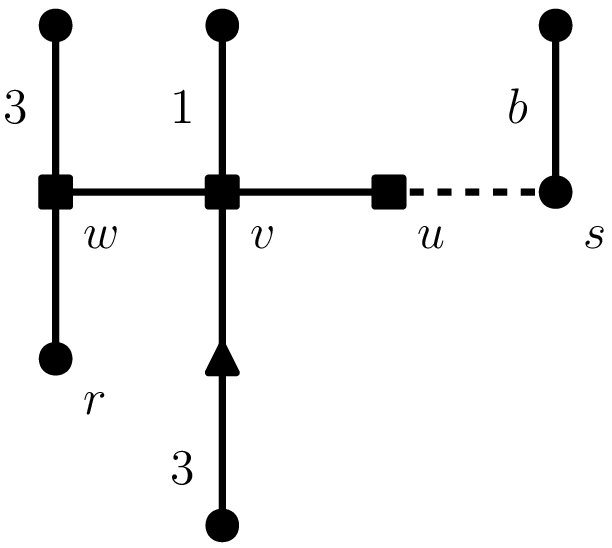}}  \hfill%
\subcaptionbox{Coloring $vu$ with 2 would lead to an unpermitted subtree. \label{fig:4lct1cse312a}} [.45\linewidth]
 {\includegraphics[scale=.5]{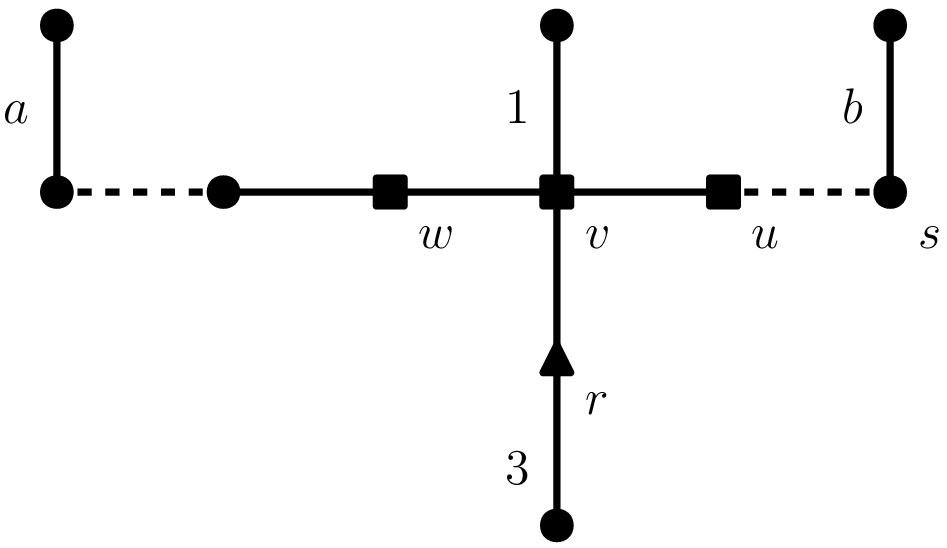}}  
	\caption{Some subtrees $F$ having exactly one colored star-edge.}
\end{figure}

If $v$ is adjacent to two 4-vertices $w$ and $u$, $F$ can be represented by the configuration in Figure~\ref{fig:4lct1cse444}. Note that in $F$, $d(r)\leq 3$ because $S_3^4$ is forbidden. In $F$, if Alice colors $vu$, the generated 4-LCT with two colored star-edges, denoted by $T_v$, may be unpermitted. Furthermore, note that in any unpermitted 4-LCT with one 4-SN and two colored star-edges, i.e., the one in Figure~\ref{fig:3124},~\ref{fig:3123} or~\ref{fig:312a}, the 4-SN is adjacent to a 4-vertex and a vertex with degree at least 3; therefore, $T_v$ is unpermitted only if $d(r)=3$. In the following, we discuss all possible situations of $F$ in which coloring $vu$ with an arbitrary color may generate an unpermitted $T_v$. For situations not mentioned below, Alice may color $vu$ with any available color to generate a Type 1 and a Type 2 subtree.

\begin{itemize}
\item If coloring $vu$ may generate the subtree in Figure~\ref{fig:3124}, $F$ can be represented by the configuration in Figure~\ref{fig:4lct1cse3124}. Alice may color $vu$ with color $d\in\{3,4\}\backslash\{b\}$ to generate a $T_v$ of Type 2.
\item If coloring $vu$ may generate the subtree in Figure~\ref{fig:3123}, $F$ can be represented by the configuration in Figure~\ref{fig:4lct1cse3123}. Note that color 3 is not available for $vu$ if and only if $u=s$ and $b=3$. Therefore, if $u=s$ and $b=3$, Alice may color $wr$ with 1 to generate a Type 1 and a Type 30 subtree; otherwise, Alice may put 3 on $vu$ to generate a Type 2 and a Type 1 subtree.
\item If coloring $vu$ may generate the subtree in Figure~\ref{fig:312a}, $F$ can be represented by the configuration in Figure~\ref{fig:4lct1cse312a}. Note that color 3 is not available for $vu$ if and only if vertex $u=s$ and $b=3$. Therefore, if $u=s$ and $b=3$, Alice may color $wv$ with 3 to generate a Type 1 and a Type 2 subtree; otherwise, Alice may put 3 on $vu$ to generate a Type 1 and a Type 2 subtree.
\end{itemize}

\item Suppose $F$ has no colored star-edges. If $v$ is not adjacent to any 4-vertex, Alice may color any star-edge to generate a Type 1 subtree and a 4-LCT with exactly one SN $v$ and colored star-edge. This 4-LCT is of Type 2 because its SN $v$ is not adjacent to any 4-vertex. Note that a 4-LCT with exactly one colored star-edge and star-node is unpermitted if and only if it is the subtree in Figure~\ref{fig:212a}, which has a star-node being adjacent to two 4-vertices. 

Similarly, if $v$ is adjacent to at least one 4-vertex, denoted by $u$, Alice may color $vu$ to generate one Type 1 subtree and one 4-LCT. This 4-LCT is of Type 2 since its star-node $v$ is adjacent to at most one 4-vertex. 

\end{itemize}

\noindent (IV) $F$ is a 4-LCT with two non-adjacent 3-SNs $v$ and $u$.\\

If $d(v)=d(u)=3$, Alice may color any edge on the star-path to generate two Type 2 subtrees.

Suppose at least one of $v$ and $u$ is of degree 4. Since $v$ and $u$ are both 3-SNs in $F$, without loss of generality, we assume that $d(v)=4$. Let $r$ be the neighbor of $v$ on the star-path. If Alice colors $vr$, she would generate a 3-LCT with star-node $v$ and the other 3-LCT with star-node $u$, denoted by $T_v$ and $T_u$, respectively. 

We first prove by contradiction that $T_u$ must be permitted, regardless of the color put on $vr$. Suppose $T_u$ is unpermitted, i.e., $T_u$ is the subtree in Figure~\ref{fig:312} or~\ref{fig:212}, in which any colored edge is incident with a 4-vertex and this 4-vertex is either the star-node or a neighbor of the star-node. Consequently, because $vr$ was colored, $r$ should be a 4-vertex in $T_u$ and $F$; moreover, $r$ is a neighbor of $u$ in $T_u$ and $F$, given that $u$ is the star-node of $T_u$. Then, in $F$, we have $d(r)=d(u)=d(s)=4$, where $s$ is a neighbor of the 3-SN $u$ and $s\neq r$, given that the star-node in an unpermitted 3-LCT is adjacent to two 4-vertices. Consequently, in $F$, the path $vrus$ is $P_4^4$, which is forbidden. 

Second, we discuss the color that Alice should put on $vr$ such that $T_v$ is also permitted. In $F$, there are exactly two colored edges which can be connected to $v$ by leaf-paths not containing $u$. We denote the color(s) on these two colored edges by $a$ and $b$. 

\bit
\item If $a$ or $b$ is available for $vr$, Alice may put $a$ or $b$ on $vr$ such that in $T_v$, a colored star-edge and a colored edge which is not a star-edge share the same color, which implies $T_v$ must not be the subtree in Figure~\ref{fig:312} nor~\ref{fig:212}. 

\item If both $a$ and $b$ are not available for $vr$, there are only two possibilities of $F$. The first one is $v$ is incident with two colored edges; then, Alice may color $vr$ with any available color to generate a $T_v$ with three colored star-edges, which is of Type 2. The second situation is $a=b$ and $v$ is incident with exactly one colored star-edge; then, Alice may also color $vr$ with any available color. After that, in $T_v$, a colored star-edge and a colored edge which is not a star-edge share the same color, which means $T_v$ is of Type 2.
\eit

\noindent (V) $F$ is a 4-LCT with two adjacent 3-SNs $v$ and $u$.\\

If Alice colors $vu$, she would generate a 3-LCT with star-node $v$ and the other 3-LCT with star-node $u$, denoted by $T_v$ and $T_u$, respectively. 

We now prove that at most one of $T_v$ and $T_u$ is unpermitted. Suppose at least one of $T_v$ and $T_u$ is unpermitted. Since $v$ and $u$ are both 3-SNs in $F$, without loss of generality, we assume that $T_v$ is unpermitted. Then, in $T_v$, the 3-SN $v$ is adjacent to two 4-vertices. Hence in $F$, $v$ is also of degree 4 and adjacent to those two 4-vertices, which means $d(u)=3$ in $F$, given that $S_3^4$ is forbidden. In addition, $d(u)=3$ implies $T_u$ must be permitted.  In the following, when one of $T_v$ and $T_u$ is unpermitted, without loss of generality, we assume $T_v$ is unpermitted.


We then discuss all situations of $F$ in which coloring $vu$ with any available color must generate an unpermitted 3-LCT, i.e., an unpermitted $T_v$. For situations not mentioned below, Alice may color $vu$ with an appropriate color to generate two 3-LCTs of Type 2.

\begin{figure}
\centering
\subcaptionbox{Coloring $vu$ must lead to the subtree in Figure~\ref{fig:312}. \label{fig:4lctadj312}} [.45\linewidth]
 {\includegraphics[scale=.5]{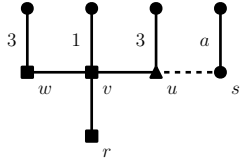}} \hfill%
\subcaptionbox{Coloring $vu$ must lead to the subtree in Figure~\ref{fig:212}. \label{fig:4lctadj212}} [.45\linewidth]
 {\includegraphics[scale=.5]{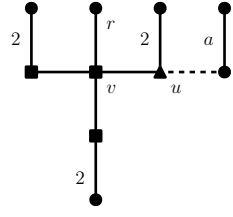}}
	\caption{Some subtrees $F$ having exactly two adjacent 3-SNs $v$ and $u$ such that coloring $vu$ must lead to an unpermitted 3-LCT}
\end{figure}

\begin{itemize}
\item If coloring $vu$ must generate the subtree in Figure~\ref{fig:312}, $F$ can be represented by the configuration in Figure~\ref{fig:4lctadj312}. Note that:
\begin{itemize}
\item Since the subtree in Figure~\ref{fig:312} has the 3-SN with two colored star-edges and $v$ would be the 3-SN in $T_v$, $v$ in $F$ should be incident with exactly one colored star-edge (without loss of generality, with color 1);
\item based on the configuration in Figure~\ref{fig:312}, $w$ should be incident with a colored edge with a color not equal to 1 (without loss of generality, with color 3);
\item $u$ should be incident with a colored edge with color 3; otherwise, Alice may color $vu$ with 3 to generate two Type 2 subtrees; and
\item if a color which is not 1 nor 3, e.g., 2, is now put on $vu$, the subtree in Figure~\ref{fig:312} will be generated.
\end{itemize}  Then, Alice's alternative move is to color $vr$ with 3 to generate a Type 1 and a Type 9 (if $s=u$) or Type 12 (if $s\neq u$) subtree.

\item If coloring $vu$ must generate the subtree in Figure~\ref{fig:212}, $F$ can be represented by the configuration in Figure~\ref{fig:4lctadj212}. Using an analysis similar to that for the last situation of $F$, we know that $u$ should be incident with a colored edge with color 2; otherwise, Alice may color $vu$ with 2 to generate two Type 2 subtrees. Moreover, if a color which is not 2, e.g., 1, is now put on $vu$, the subtree in Figure~\ref{fig:212} will be generated. Alice's alternative move is to color $vr$ with 2 to generate a Type 1 and a Type 37 subtree.
\end{itemize}

\noindent (VI) $F$ is a 5-LCT with a 4-SN $v$ and 3-SN $u$ which are not adjacent. \\

\begin{figure}
\centering
\subcaptionbox{The 4-SN $v$ and 3-SN $u$. \label{fig:5lctnonadj}} [.45\linewidth]
 {\includegraphics[scale=.5]{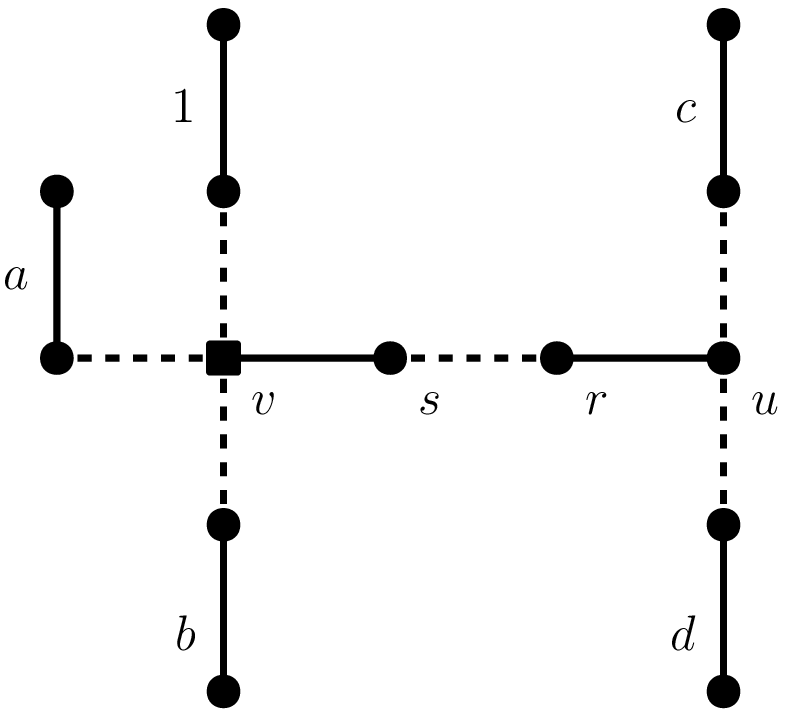}} \hfill%
\subcaptionbox{Coloring $vs$ with 4 would lead to an unpermitted subtree. \label{fig:5lctnonadj1234}} [.45\linewidth]
 {\includegraphics[scale=.5]{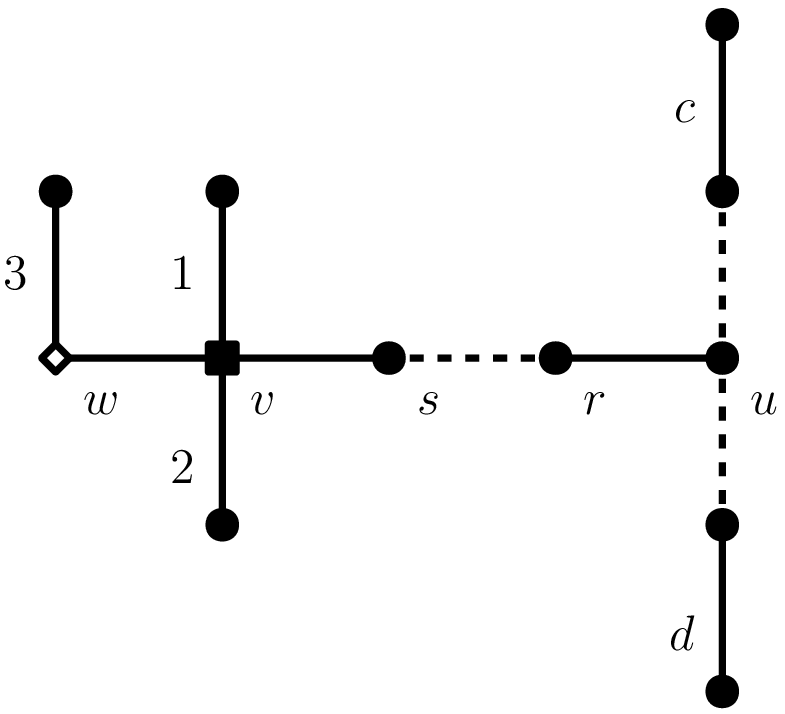}}
 
\subcaptionbox{Coloring $vs$ with 3 would lead to an unpermitted subtree. \label{fig:5lctnonadj3123}} [.45\linewidth]
 {\includegraphics[scale=.5]{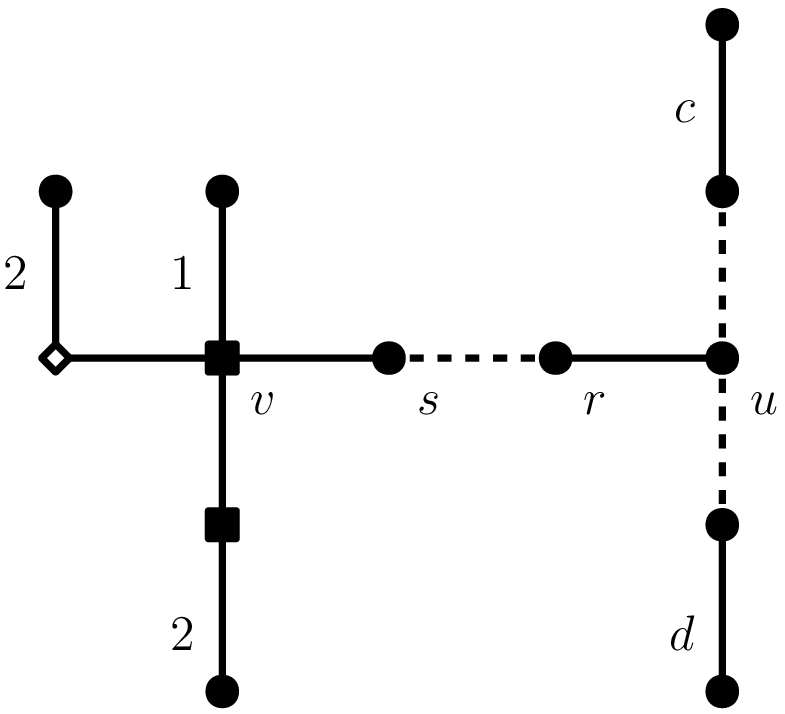}}\hfill%
\subcaptionbox{Coloring $vs$ with 4 would lead to an unpermitted subtree. \label{fig:5lctnonadj3124}} [.45\linewidth]
 {\includegraphics[scale=.5]{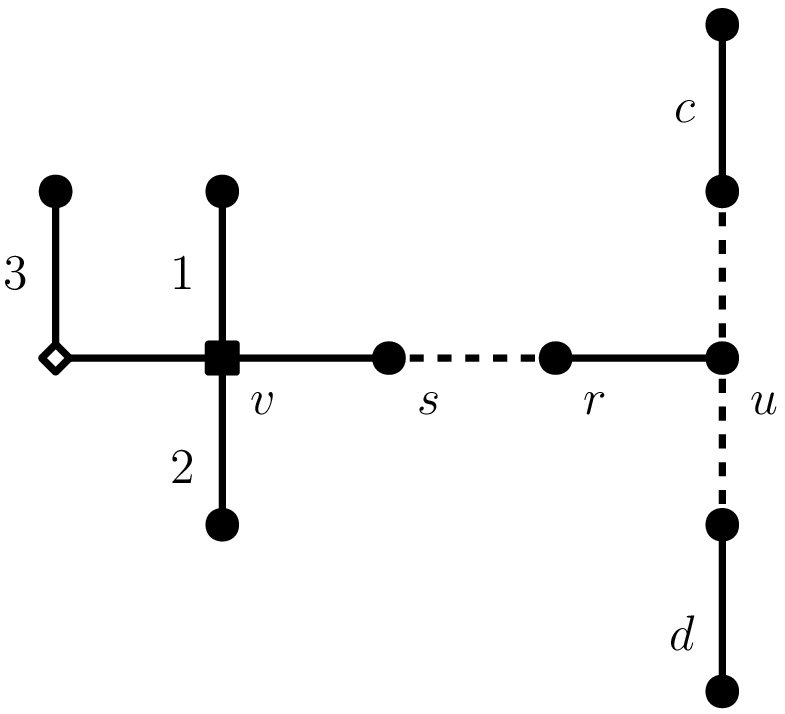}}
 
\subcaptionbox{Coloring $vs$ with 3 would lead to an unpermitted subtree. \label{fig:5lctnonadj312a}} [.45\linewidth]
 {\includegraphics[scale=.5]{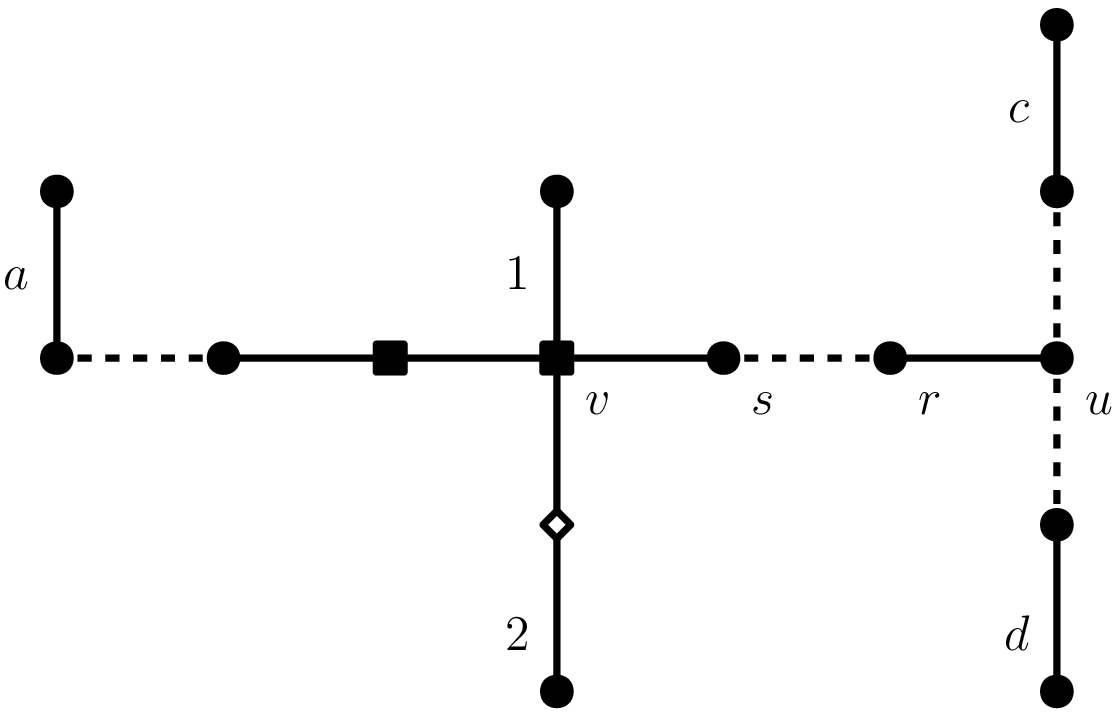}}\hfill%
\subcaptionbox{Coloring $vs$ with 2 would lead to an unpermitted subtree. \label{fig:5lctnonadj212a}} [.45\linewidth]
 {\includegraphics[scale=.5]{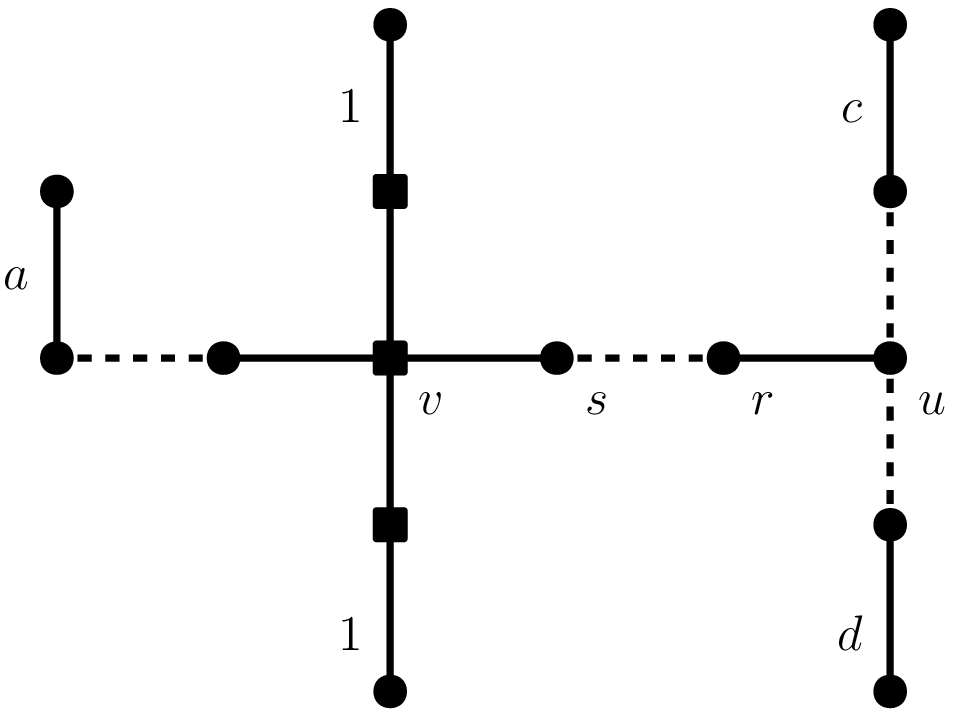}}
	\caption{Some subtrees $F$ having exactly two non-adjacent star-nodes}
\end{figure}

$F$ in this case can be represented by the configuration in Figure ~\ref{fig:5lctnonadj}. If $vs$ is colored, $F$ is split into one 4-LCT with 4-SN $v$ (denoted by $T_v$) and one 3-LCT with 3-SN $u$ (denoted by $T_u$). We will first prove that no matter what color is put on $vs$, $T_u$ is of Type 2, i.e., not the subtree in Figure~\ref{fig:312} nor~\ref{fig:212}. After that, we will describe, under different situations of $F$, what color Alice should put on $vs$ such that $T_v$ is also of Type 2. These mean that Alice can always decompose $F$  into two Type 2 subtrees by coloring $vs$.

We now prove by contradiction that $T_u$ is of Type 2. Suppose $T_u$ is the subtree in Figure~\ref{fig:312} or~\ref{fig:212}. Then, in $F$, $s=r$ and $d(s)=4$; also, in $F$, $d(u)=4$ and $u$ is adjacent to another 4-vertex $x$ which is not the vertex $s$. Hence the path $vsux$ is $P_4^4$, which is forbidden.

Then, we discuss all possible situations of $F$ in which coloring $vs$ with an arbitrary color may generate an unpermitted $T_v$, i.e., the subtree in Figure~\ref{fig:1234},~\ref{fig:3123},~\ref{fig:3124},~\ref{fig:312a} or~\ref{fig:212a}. We will show that Alice can always find a suitable color $c$ for coloring $vs$ such that the generated $T_v$ is of Type 2. Because in the subtrees in Figure~\ref{fig:1234}-~\ref{fig:312a}, any colored star-edge does not share the same color with any colored edge which is incident with a neighbor of the star-node $v$, Alice's requirement on $c$ is it should be a color on a colored edge which is incident with a neighbor of $v$. For situations not mentioned below, Alice may color $vs$ with any available color to generate two Type 2 subtrees. 

\begin{itemize}
\item If coloring $vs$ may generate the subtree in Figure~\ref{fig:1234}, $F$ can be represented by the subtree in Figure~\ref{fig:5lctnonadj1234}. Then, Alice may put 3 on $vs$.
\item If coloring $vs$ may generate the subtree in Figure~\ref{fig:3123}, $F$ can be represented by the subtree in Figure~\ref{fig:5lctnonadj3123}. Then, Alice may put 2 on $vs$.
\item If coloring $vs$ may generate the subtree in Figure~\ref{fig:3124}, $F$ can be represented by the subtree in Figure~\ref{fig:5lctnonadj3124}. Then, Alice may put 2 or 3 on $vs$.
\item If coloring $vs$ may generate the subtree in Figure~\ref{fig:312a}, $F$ can be represented by the subtree in Figure~\ref{fig:5lctnonadj312a}. Then, Alice may put 2 on $vs$.
\item If coloring $vs$ may generate the subtree in Figure~\ref{fig:212a}, $F$ can be represented by the subtree in Figure~\ref{fig:5lctnonadj212a}. Then, Alice may put 1 on $vs$.
\end{itemize}

\noindent (VII) $F$ is a 5-LCT with a 4-SN $v$ and 3-SN $u$ which are adjacent. \\

We prove by contradiction that there is always an available color for $vu$, that is, $vu$ should not be surrounded by edges with five colors on them. Suppose it is not. Before Bob's last move, at least four colors were on the edges surrounding $vu$, which implies the subtree was the unpermitted one in Figure~\ref{fig:1234}. This is impossible since an unpermitted $F$ must be generated from a permitted subtree by Bob.

If Alice colors $vu$, she will generate one 4-LCT with 4-SN $v$ (denoted by $T_v$) and one 3-LCT with 3-SN $u$ (denoted by $T_u$). So, Alice may color $vu$ if only if $T_v$ and $T_u$ will be both permitted. We first show that $T_u$ must be of Type 2. If $d(u)=3$, $T_u$ is of Type 2. If $d(u)=4$, since $S_3^4$ is forbidden and $d(v)=4$, we conclude that the 3-SN $u$ is adjacent to at most one 4-vertex in $T_u$, which means $T_u$ is of Type 2.

$T_v$ is of Type 2 unless it is the subtree in Figure~\ref{fig:1234},~\ref{fig:3123},~\ref{fig:3124},~\ref{fig:312a} or~\ref{fig:212a}. In the following, we will discuss all possible configurations of $F$ in which coloring $vs$ must lead to an unpermitted $T_v$, and describe the corresponding alternative moves of Alice. For situations not mentioned below, Alice may color $vu$ with an appropriate color to generate two Type 2 subtrees.

\begin{itemize}

\begin{figure}
\centering
\subcaptionbox{Coloring $vu$ with 4 would lead to the subtree in Figure~\ref{fig:1234}. Vertex $s$ exists if $d(u)=4$. \label{fig:5lctadj1234}} [.45\linewidth]
 {\includegraphics[scale=.5]{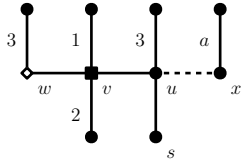}}\hfill%
\subcaptionbox{Coloring $vu$ with 3 would lead to the subtree in Figure~\ref{fig:3123}. \label{fig:5lctadj3123}} [.45\linewidth]
 {\includegraphics[scale=.5]{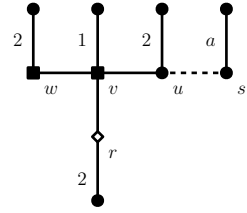}}

\subcaptionbox{Coloring $vu$ with 4 would lead to the subtree in Figure~\ref{fig:3124}. \label{fig:5lctadj3124}} [.45\linewidth]
 {\includegraphics[scale=.5]{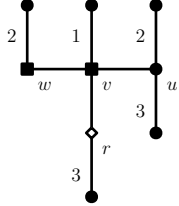}}\hfill%
 \subcaptionbox{Coloring $vu$ with 3 would lead to the subtree in Figure~\ref{fig:312a}. \label{fig:5lctadj312a}} [.45\linewidth]
 {\includegraphics[scale=.5]{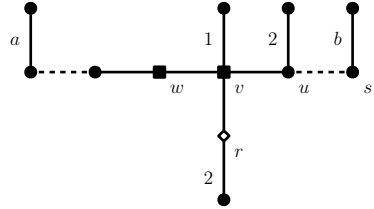}}

\subcaptionbox{Coloring $vu$ with 2 would lead to the subtree in Figure~\ref{fig:212a}. Since $S_3^4$ is forbidden, $d(u)=3$ and $d(t)\leq 3$. \label{fig:5lctadj212a}} [.45\linewidth]
 {\includegraphics[scale=.5]{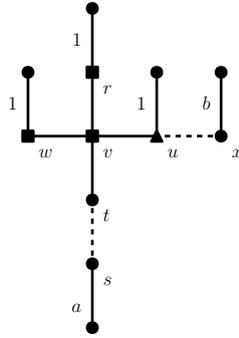}}
	\caption{Some subtrees $F$ having exactly two adjacent star-nodes, which are the 4-SN $v$ and 3-SN $u$.}
\end{figure}

\item If coloring $vu$ must lead to the subtree in Figure~\ref{fig:3124}, $F$ can be represented by the subtree in Figure~\ref{fig:5lctadj3124}. In $F$, note that $u$ should be incident with edges colored with 2 and 3 so that coloring $vu$ with 4 or 5 would lead to the subtree in Figure~\ref{fig:3124}. We discuss this $F$ with two cases: (1) $d(r)=3$ and (2) $d(r)=4$.

\begin{enumerate}
\item Given that $r$ is a 3-vertex, Alice may color $wv$ with 3 to generate one Type 1 and one Type 4 (when $d(u)=3$) or Type 6 (when $d(u)=4$) subtree.
\item Given that $r$ is a 4-vertex and $S_3^4$ is forbidden, we conclude that $d(u)=3$. Alice may also color $wv$ with 3 to generate one Type 1 and one Type 9 subtree.
\end{enumerate}

\item If coloring $vu$ must lead to the subtree in Figure~\ref{fig:312a}, $F$ can be represented by the subtree in Figure~\ref{fig:5lctadj312a}. In $F$, note that $u$ should be incident with an edge with color 2 so that coloring $vu$ with 3/4/5 would lead to the subtree in Figure~\ref{fig:3124}. Because $S_3^4$ is forbidden, at most one of $r$ and $u$ is of degree 4. 

Suppose $u=s$. Alice may color $wv$ with 2 to generate a Type 1 subtree, and a Type 9 (if $d(u)=3$) or Type 6 (if $d(u)=4$) subtree.

Suppose $u\neq s$. Alice may also color $wv$ with 2 to generate a Type 1 subtree, and a Type 38 (if $d(r)=3$) or Type 12 (if $d(r)=4$) subtree, given that $S_3^4$ and $P_4^4$ are forbidden.

\item If coloring $vu$ must lead to the subtree in Figure~\ref{fig:212a}, $F$ can be represented by the subtree in Figure~\ref{fig:5lctadj212a}. In $F$, note that each of $u,w$ and $r$ is incident with an colored edge with the same color, i.e., color 1. Also, because $S_3^4$ is forbidden, we have $d(u)=3$ and $d(t)\leq 3$.

If $t$ is not incident with an edge with color 1, Alice may color $vt$ with 1 to generate a Type 1 and a Type 37 subtree.

Suppose $t$ is incident with an edge with color 1, i.e., $t=s$ and $a=1$. Then, $F$ is now of Type 21 if $u=x$. If $u\neq x$, we let $y$ be the neighbor of $u$ on the path connecting $u$ and $x$. Alice may  color $uy$ with any available color to generate a Type 21 and a Type 1 subtree.

\item If coloring $vu$ must lead to the subtree in Figure~\ref{fig:1234}, $F$ can be represented by the subtree in Figure~\ref{fig:5lctadj1234}. In $F$, note that $u$ must be incident with an edge colored with 3 so that only colors 4 and 5 are available for $vu$. We now describe alternative moves of Alice for two cases: (1) $u$ is incident with two colored edges (i.e., $x=u$) and (2) $u$ is incident with one colored edge ($x\neq u$).

\ben
\item When $x=u$, we consider two situations: (i) $a\neq1,2$ and (ii) $a=1$ or 2. 

(i) Without loss of generality, assume $a=4$. Note that $d(w)\geq 3$, and a color of the colored star-edges incident with $u$ was just added by Bob in his last move on a Type 2 subtree. Since the trees in Figure~\ref{fig:3123} and~\ref{fig:3124} are unpermitted, we have $d(u)\leq 3$, which implies $d(u)=3$. Alice may put 4 on $wv$ to generate one Type 1 and one Type 3 subtree. 

(ii) If $a=1$ or 2, without loss of generality, assume $a=1$. For the case that $d(u)=3$, Alice may put 4 on $wv$ to generate one Type 1 and one Type 3 subtree. If $d(u)=4$, Alice may put 2 on $us$ to generate one Type 1 and one Type 40 subtree. 

\item When $x\neq u$ and $d(u)=3$, Alice may put 1 or 2 (at least one of them is available) on the star-edge incident with $u$ on the leaf-path connecting $u$ and $x$ to generate one Type 1 and Type 6 subtree. We will show that $d(u)\neq 4$, using the fact that either the color 3 of the edge incident with $u$ or the color $a$ was just added by Bob in his last move on a permitted Type 2 subtree. 

Suppose the color 3 of the edge incident with $u$ was just added. Since $d(w)\geq 3$, $u$ cannot be a 4-vertex; otherwise, the subtree after the color 3 is removed will be the one in Figure~\ref{fig:312a}, which contradicts that any unpermitted $F$ was generated from a permitted subtree. 

Similarly, if the color $a$ was just added, $u$ cannot be a 4-vertex; otherwise, the subtree after the color $a$ is removed will be the one in Figure~\ref{fig:3123}.
\een

\item If coloring $vu$ must lead to the subtree in Figure~\ref{fig:3123}, $F$ can be represented by the subtree in Figure~\ref{fig:5lctadj3123}. In $F$, note that $u,w,r$ should have the same color (i.e., color 2) on the edges incident with them so that coloring of $vu$ with 3, 4 or 5 would lead to the subtree in Figure~\ref{fig:3123}. Note that  and either the color 2 of the colored edge incident with $u$ or the color $a$ of the colored edge incident with $s$ was just added by Bob. Because the subtree in Figure~\ref{fig:212a} is unpermitted, we conclude that at most one of $w$ and $r$ is of degree 4, which implies $d(r)=3$. 

In the following, we will introduce alternative strategies of Alice for two cases: (1) $u=s$ and (2) $u\neq s$.

\begin{enumerate}
\item If $u=s$, then $s$ is now incident with two edges with colors 2 and $a$; hence we may assume $a=1$ or 3. Alice may color $wv$ with 3 to generate a Type 1 and a Type 4 (if $d(u)=3$) or Type 6 (if $d(u)=4$).

\item Suppose $u\neq s$. Let $t$ be the neighbor of $u$ which is on the path connecting $u$ and $s$. We consider two situations: (i) $d(u)=3$ and (ii) $d(u)=4$. 

(i) Suppose $d(u)=3$. If $t$ is not incident with an edge colored with 1, Alice may color $ut$ with 1 to generate a Type 8 and a Type 1 subtree. If $t$ is incident with an edge colored with 1, i.e., $t=s$ and $a=1$, Alice may color $wv$ with 3 to generate a Type 1 and a Type 4 (if $d(t)=2$) or 17 (if $d(t)\geq 3$) subtree.

(ii) When $d(u)=4$, we let $x$ be the neighbor of $u$ which is of the uncolored $u$-branch. If $t$ is not incident with an edge colored with 1, Alice may color $ut$ with 1 to generate a Type 30 and a Type 1 subtree. If $t$ is incident with an edge colored with 1, i.e., $t=s$ and $a=1$, Alice may color $ux$ with 1 to generate a Type 34 and a Type 1 subtree.
\end{enumerate}

%
%
%

\end{itemize} 

\section{Proof of Lemma~\ref{lem:inductive}: Alice's strategies for Types 3 to 40 subtrees}
\medskip
We will give, one by one, Alice's strategies for handling (i) subtrees of each of Types 3 to 40; and (ii) any unpermitted subtree generated by Bob's move on subtrees of Types 3 to 40, respectively. Alice may need to use an appropriate color to generate the desired types of subtrees.
\medskip

\noindent\tb{Type 3 subtrees:} {\it [Recall the definition: A subtree with the subgraph induced by all its uncolored edges being a path $P$ of length $m$, where $1\leq m\leq 3$. Moreover, for $m>1$, at least two colors are available for each uncolored edge; the only uncolored edge has at least one available color for $m=1$.]}

(i) When Alice is going to act on a Type 3 subtree, if $m=3$, she may color the middle edge to generate two Type 3 subtrees of length 1. If $m<3$, she may color any edge to generate a Type 0 and a Type 3 subtree (for $m=2$), or two Type 0 subtrees (for $m=1$).

(ii) If $m=3$ and Bob has made an unpermitted subtree by acting on an end edge of $P$, Alice may color the middle edge of $P$ to generate a Type 0 subtree and a Type 3 subtree with $m=1$.\\

\noindent\tb{Type 4 subtrees:} {\it [Recall the definition: A subtree with the subgraph induced by all its uncolored edges being the union of a path $P=v_0v_1...v_m$, where $3\leq m\leq 4$, and a tree $T_s$ (maybe trivial). Also, this subtree has to satisfy the two requirements that (1) $v_m$ is the only common vertex of $P$ and $T_s$, and (2) all vertices of $T_s$ are not incident with any colored edge. Moreover, the path $P$ has to satisfy the condition that $v_{m-1}v_m$ has at least four available colors and each of the other uncolored edges of $P$ has at least two available colors.]}

(i) We observe that in all Type 4 subtrees, $d(v_{m-1})\leq 3$ and $v_{m-1}$ is incident with at most one colored edge. Alice may put an appropriate color on $v_{m-1}v_{m}$ such that at least two colors are still available for $v_{m-2}v_{m-1}$ to generate a Type 1 and a Type 3 subtree. 

(ii) After Bob has acted on a Type 4 subtree, Alice may respond as follows:

\begin{itemize}
\item If Bob has dyed $v_0v_1$ (resp. $v_2v_3$), Alice may dye $v_2v_3$ (resp. $v_0v_1$) to make a Type 1 subtree containing $v_2v_3$ (resp. a Type 0 subtree containing $v_0v_1$), and a Type 3 subtree with only one uncolored edge $v_1v_2$. 

\item If he has dyed $v_3v_4$ (only for $m=4$), Alice may dye $v_1v_2$ to make a Type 3 subtree with only one uncolored edge $v_0v_1$ and the other Type 3 subtree with only one uncolored edge $v_2v_3$.

\item If he has dyed $v_1v_2$, Alice may dye $v_{m-1}v_m$ to make a Type 0 subtree containing $v_1v_2$ and $v_2v_3$ (when $m=3$) or a Type 3 with only one uncolored edge $v_2v_3$ (when $m=4$), and a Type 1 with only one colored edge $v_{m-1}v_m$.

\item If he has dyed anyone of other edges, Alice may dye $v_{m-1}v_m$ to make a Type 3 subtree with only two uncolored edges (when $m=3$) or only three (when $m=4$), and a Type 1 subtree with only two colored edges.
\end{itemize}

\noindent\tb{Types 5, 19 and 20 subtrees:}

We propose a unified strategy of Alice for any subtree $H$ of Type 5, 19 or 20. Note that since $P_4^4$ is forbidden, $d(s)\leq 3$ for any Type 20 subtree. Also, at least colors 4 and 5 are available for $wv, vu$ and $us$. Therefore, if $d(s)=1$ (resp. $d(s)=2$), $H$ is of Type 3 (resp. Type 4 with $m=4$) and Alice can use the corresponding strategies to handle $H$. So, we assume $d(s)=3$ in the following.

(i) Alice can color $vu$ to generate one Type 3 and one Type 2 subtree.

(ii) After Bob has colored an edge of $H$, Alice may respond as follows:

\begin{itemize}
\item If Bob colored any edge incident with $s$, except $us$, he generated a Type 4 and a Type 1 subtree. So, all these generated subtrees are permitted and Alice can use the corresponding strategies to handle these subtrees.

\item If Bob colored $us$, Alice may color $vu$ to generate one Type 3 and one Type 0 subtree.  

\item If Bob colored anyone of other edges, Alice may color $us$ with an available color in $\{1,2,3\}$ to generate a Type 3 subtree and a Type 1 subtree with at most two colored edges.
\end{itemize}

\noindent\tb{Type 6 subtrees:}

(i) Alice may color $us$ with any available color in $\{1,2,3\}$ to generate one Type 19 and one Type 1 subtree.

(ii) If Bob has colored an edge, Alice may respond as follows:

\begin{itemize}
\item If Bob colored $us$ (resp. $wv$) with 4, Alice may put 4 on $wv$ (resp. $us$) to generate a Type 2 and Type 3 subtree.
\item If he colored $vu$ (resp. $rw$) with $a$, she may color $rw$ (resp. $vu$) with $a$ if available or other colors if not to generate a Type 1 and a Type 3 subtree. 
\item If he colored an edge incident with $r$, except $rw$, she may color $us$ with any available color in $\{1,2,3\}$ to generate a Type 3/4 (since $d(r)\leq 3$) and a Type 1 subtree. If Bob colored another edge of an uncolored $r$-branch, Alice may color $rw$ with any available color in $\{1,2,3\}$ to generate one Type 1 and one Type 4 subtree.
\item If he colored an edge of an uncolored $s$-branch, she may color $us$ with any available color in $\{1,2,3\}$ to generate one Type 19 and one Type 1 subtree.
\end{itemize}



\noindent\tb{Types 7 and 39 subtrees:}

Note that $d(w)=4$ in Type 39 subtrees, which implies vertices $r,s$ and $t$ are not of degree 4, given that $S_3^4$ and $P_4^4$ are forbidden. Therefore, in both Type 7 and 39 subtrees, vertices $r,s$ and $t$ are not of degree 4. We will provide a unified strategy of Alice for any subtree $H$ which is of Type 7 or 39. Note that vertices $t$ and $s$ are identical by symmetry.

(i) Alice may color $vu$ with 3 to generate a Type 1 and a Type 9 (if $H$ is of Type 7) or 40 (if $H$ is of Type 39) subtree.

(ii) If Bob has colored an edge, Alice may respond as follows.
\begin{itemize}
\item If he colored $wv$ with 4, she may color $us$ with 4 to generate a Type 1 and a Type 4 (if $d(r)=2$) or Type 6 (if $d(r)=3$) subtree.
\item If he colored $vr$ with 3 (resp. 4), she may color $vu$ with 4 (resp. 3) to generate a Type 3 and a Type 1 subtree.
\item Suppose $d(r)=3$. If he colored $rx$ with 3, she may color $vu$ with 3 to generate a Type 3 and a Type 1 subtree. 
\item Suppose $d(r)=3$. If he colored $rx$ (resp. an edge in an uncolored $s$-branch) with $a\neq 3$ , she may color $vr$ with 3 to generate a Type 9 (when $H$ is of Type 7) or Type 40 (when $H$ is of Type 39) subtree and a Type 0 (resp. Type 2) subtree.
\item If he colored $vu$ with 4, she may color $vr$ with 3 to generate a Type 3 and a Type 1 subtree.
\item If he colored $us$ with 3, she may color $vr$ with 3 to generate a Type 1 and a Type 19 (if $H$ is of Type 7) or Type 20 (if $H$ is of Type 39) subtree.
\item If he colored $us$ with color $a\neq3$ or an edge in an uncolored $s$-branch, she may color $vu$ with 3 to generate a Type 9 (if $H$ is of Type 7) or Type 40 (if $H$ is of Type 39) subtree and a 3-LCT with 3-SN $u$, denoted by $T_u$. In $T_u$, because $d(s),d(t)\leq 3$, $u$ is not adjacent to two 4-vertices, which means $T_u$ is of Type 2.
\end{itemize}

\noindent\tb{Type 8 subtrees:}

(i) Alice may color $vu$ with 3 to generate a Type 1 and a Type 9 subtree.

(ii) If Bob has colored an edge, Alice may respond as follows.
\begin{itemize}
\item If he colored $wv$ with 3, she may color $rs$ with 1 to generate a Type 9 and a Type 1 subtree.
\item If he colored $vr$ with 3, he generated a Type 1 and a Type 9 subtree.
\item If he colored an edge in an uncolored $r$-branch, she may color $vr$ to generate a Type 9 and a Type 0/2 subtree.
\item If he colored an edge in an uncolored $u$-branch, she may color $vu$ to generate a Type 9 and a Type 2 (since $S_3^4$ is forbidden) subtree.
\end{itemize}

\noindent\tb{Types 9 and 40 subtrees:}

We will propose a unified strategy for both Type 9 and 40 subtrees.

(i) Alice may color $vu$ with any available color to generate one Type 3 and one Type 1 subtree. 

(ii) If Bob has colored an edge, Alice may respond as follows:
\begin{itemize}
\item If he colored $wv$, Alice may color $vu$ to generate a Type 0 and a Type 1 subtree. 
\item If he colored an edge of a $u$-branch not containing $v$, Alice may color $vu$ with any available color to generate a Type 3 subtree and a 2-/3-LCT containing $u$, denoted by $T_u$. $T_u$ is of Type 2 if it has two colored edges. Suppose $T_u$ is a 3-LCT. Because $S_3^4$ and $P_4^4$ are forbidden, in $T_u$, $u$ is neither the middle nor the end vertex of a $P_3^4$; therefore, $T_u$ is neither the subtree in Figure~\ref{fig:212} nor~\ref{fig:312}, which implies $T_u$ is of Type 2.
\end{itemize}

\noindent\tb{Type 10 subtrees:}

(i) Alice may color $vu$ or $vr$ with any available color to generate a Type 9 and a Type 1 subtree.

(ii) If Bob has colored an edge, Alice may respond as follows:
\bit
\item If he colored $wv$, she may color $vu$ with any available color to generate a 4-LCT of Type 2 and a Type 1 subtree.
\item If he colored an edge of a $u$-branch not containing $v$, she may color $vu$ with any available color to generate a Type 9 subtree and a 2-/3-LCT containing $u$, denoted by $T_u$. $T_u$ is of Type 2 if it has two colored edges. Suppose $T_u$ is a 3-LCT. Because $S_3^4$ and $P_4^4$ are forbidden, in $T_u$, $u$ is neither the middle nor an end vertex of a $P_3^4$; therefore, $T_u$ is neither the subtree in Figure~\ref{fig:212} nor~\ref{fig:312}, which implies $T_u$ is of Type 2. 

The case that Bob colored an edge of a $r$-branch not containing $v$ can be analyzed similarly.
\eit

\noindent\tb{Types 11, 14, 15, 16 and 18 subtrees:}

We provide a unified strategy for all the five types of subtrees. Denote the subtree that Bob is acting on by $Q$. Note that $tw$ and $ur$ are identical by symmetry when $Q$ is of Type 15 or 18; and $wx$ and $us$ are identical by symmetry when $Q$ is of Type 18.

(i) When $Q$ is of Type 11, 14, 15, 16 or 18, Alice may color $ur$ with any available color to generate a Type 1 subtree, and a Type 3, 11, 11, 15 or 16 subtree, respectively.

(ii) If Bob has colored an edge of $Q$, Alice may respond as follows:
\bit
\item Suppose Bob colored $vu$. Then, Alice may color $ur$ with any available color to generate a Type 1 and a Type 0/2 subtree. 

Similarly, suppose Bob colored $wv$. If not all subtrees are permitted, then $tw$ exists and Alice may color $tw$ with any available color to generate a Type 1 and a Type 0/2 subtree.

\item If he colored $tw$, $us$ or $wx$, all subtrees are now permitted.

\item Suppose he colored an edge of a $r$-branch not containing $u$. She may color $ur$ to generate a 2-/3-LCT (denoted by $T_r$) and a Type 3, 11, 11, 15 or 16 subtree. $T_r$ is of Type 1 if it is a 2-LCT. Suppose $T_r$ is a 3-LCT. Given that $S_3^4$ and $P_4^4$ are forbidden and $d(u)=4$, $r$ is neither the middle nor the end vertex of a $P_3^4$ in $T_r$, which implies $T_r$ is of Type 2.

Similar analysis can be done when Bob colored an edge of a $s$-branch not containing $u$, a $t$-branch not containing $w$, or a $x$-branch not containing $w$.
\eit

\noindent\tb{Type 12 subtrees:}

(i) Alice may color $vu$ to generate two Type 2 subtrees.

(ii) If Bob has colored an edge, Alice may respond as follows:
\bit
\item If he colored an edge of an uncolored $w$-branch, she may color $wv$ to generate a 3-LCT with star-node $w$ and a Type 13 subtree. Because $S_3^4$ is forbidden, in the 3-LCT, $w$ is adjacent to at most one 4-vertex, which implies this 3-LCT is of Type 2. 
\item If he colored $wv$ or $ut$, Alice may color $vu$ to generate a Type 0 and a Type 2 subtree.
\item Suppose he colored an edge of a $t$-branch not containing $u$. Then, a Type 12 subtree was split into a Type 1 subtree and the other subtree, denoted by $H$. We will describe possible moves of Alice on $H$ in the following. If Alice colors $ut$, she will generate a Type 9 subtree and a 2-/3-LCT, denoted by $T_t$. $T_t$ is of Type 1 if it is a 2-LCT; and of Type 2 if and only if it is a 3-LCT but not the subtree in Figure~\ref{fig:312} nor~\ref{fig:212}. Therefore, Alice may color $ut$ with a suitable color if and only if $T_t$ will be permitted. 

\begin{figure}
\centering
\subcaptionbox{Coloring $ut$ with 1 or 3 would lead to an unpermitted subtree. \label{fig:type12out212}} [.45\linewidth]
 {\includegraphics[scale=.5]{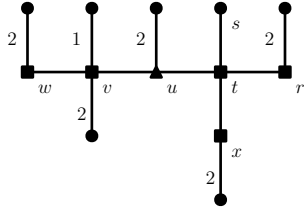}} \hfill%
\subcaptionbox{Coloring $ut$ with 1 or 3 would lead to an unpermitted subtree. \label{fig:type12out312}} [.45\linewidth]
 {\includegraphics[scale=.5]{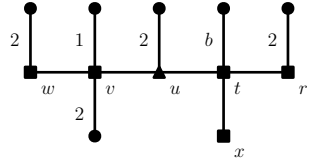}}
	\caption{Situations of Type 12 subtrees in which coloring $ut$ with any available color must lead to an unpermitted subtree.}
\end{figure}

\bit
\item If coloring $ut$ must lead to the subtree in Figure~\ref{fig:212}, $H$ can be represented by the subtree in Figure~\ref{fig:type12out212}. Note that in $H$, $r$ and $x$ should be incident with an edge colored with 2 so that Alice cannot put 2 on $ut$. Then, Alice may color $ts$ with 2 to generate a Type 1 and a Type 16 subtree.
\item If coloring $ut$ must lead to the subtree in Figure~\ref{fig:312}, $H$ can be represented by the subtree in Figure~\ref{fig:type12out312}. Note that in $H$, $r$ should be incident with an edge colored with 2 so that Alice cannot put 2 on $ut$; also, $t$ should be incident with an edge colored with $b\neq 2$. Then, Alice may color $tx$ with 2 to generate a Type 1 and a Type 15 subtree.
\eit
\eit

\noindent\tb{Type 13 subtrees:}

The strategy of Alice for this type is very similar to that for Type 12 subtrees. 

(i) Alice may color $vu$ to generate a Type 0 and a Type 2 subtree.

(ii) If Bob has colored an edge, Alice may respond as follows:
\begin{itemize}
\item If he colored $ut$, she may color $vu$ to generate two Type 0 subtrees.
\item Suppose he colored an edge of a $t$-branch not containing $u$. Then, a subtree of Type 13 will be split into a Type 1 subtree and the other subtree, denoted by $H$. We will describe possible moves of Alice on $H$ in the following. If Alice colors $ut$, she will generate a Type 3 subtree and a 2-/3-LCT with star-node $t$, denoted by $T_t$. $T_t$ is of Type 1 if it is a 2-LCT; and of Type 2 if and only if it is a 3-LCT but not the subtree in Figure~\ref{fig:212} nor~\ref{fig:312}. Therefore, Alice may color $ut$ with a suitable color if and only if $T_t$ is permitted. Now, we discuss situations in which coloring $ut$ must lead to an unpermitted $T_t$.

\begin{figure}
\centering
\subcaptionbox{Coloring $ut$ with 1/3/4 would lead to an unpermitted subtree. \label{fig:type13out212}} [.45\linewidth]
 {\includegraphics[scale=.5]{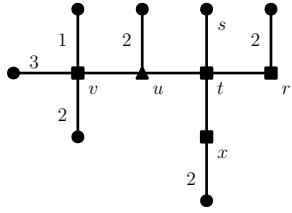}} \hfill%
\subcaptionbox{Coloring $ut$ with 1/3/4 would lead to an unpermitted subtree. \label{fig:type13out312}} [.45\linewidth]
 {\includegraphics[scale=.5]{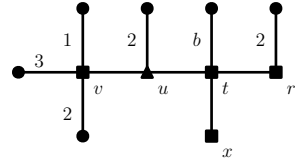}}
	\caption{Situations of Type 13 subtrees in which coloring $ut$ with any available color must lead to an unpermitted subtree.}
\end{figure}

If coloring $ut$ must lead to the subtree in Figure~\ref{fig:212} (resp. Figure~\ref{fig:312}), $H$ can be represented by the subtree in Figure~\ref{fig:type13out212} (resp. Figure~\ref{fig:type13out312}). So, Alice may color $ts$ with 2 to generate a Type 1 and a Type 14 (resp. Type 11) subtree.

\end{itemize}

\noindent\tb{Type 17 subtrees:}

Note that $rs$ and $ry$ are identical by symmetry.

(i) Alice may color $ur$ with 3 to generate a Type 4 and a Type 1 subtree.

(ii) If Bob has colored an edge, Alice may respond as follows.
\begin{itemize}
\item If he colored an edge in an uncolored $t$-branch, she may colored $tw$ with an available color in  $\{1,3\}$ to generate a Type 1 and a Type 22 subtree.
\item If he colored $tw$ with 1/3/4 (resp. $vu$ with 4), she may color $vu$ with 4 (resp. $tw$ with 1) to generate a Type 3 and a Type 2 subtree.
\item If he colored $wv$ (resp. $ur$) with 4, she may color $ur$ (resp. $wv$) with 4 to generate a Type 1 and a Type 3 subtree.
\item Suppose $d(r)=3$. If he colored $rs$, he generated a Type 4 and a Type 1 subtree. If he colored an edge in an uncolored $s$-branch, she may color $rs$ with an available color in  $\{2,3\}$ to generate a Type 4 and a Type 1 subtree.
\item Suppose $d(r)=4$. If he colored $rs$ with 2 (resp. 3 and 4), she may color $ry$ with 3 (resp. 2 and 2) to generate a Type 4 and a Type 1 subtree.
\item Suppose $d(r)=4$. If he colored an edge in an uncolored $s$-branch. Then, a Type 17 subtree was split into a Type 1 subtree and the other subtree $H$. If Alice now colors $ur$ with 3 in $H$, she will generate a Type 4 subtree, and a 3-LCT $T_r$. $T_r$ is of Type 2 if and only if it is not the subtree in Figure~\ref{fig:312} nor~\ref{fig:212}. Because $T_r$ has exactly two colored star-edges, $T_r$ must not be the one in Figure~\ref{fig:212}. Therefore, Alice may color $ur$ with 3 if and only if it will not lead to the subtree in Figure~\ref{fig:312}. Note that $T_r$ would be the subtree in Figure~\ref{fig:312} if and only if $d(s)=d(y)=4$ and $sz$ was colored with 2 or 4, since the colored star-edges already have colors 1 and 3. So, if $sz$ was colored with 2 (resp. 4), an alternative move of  Alice is to color $ry$ (resp. $rs$) with 2 to generate a Type 23 (resp. Type 24) and a Type 1 subtree.
\end{itemize}

\noindent\tb{Type 21 subtrees:}

Note that $vu$ and $vs$ are identical by symmetry.

(i) Alice may color $vu$ with 2 to generate a Type 7 and a Type 1 subtree.

(ii) If Bob has colored an edge, Alice may respond as follows.
\begin{itemize}
\item If he colored $wv$ with 3, she may color $ut$ with 3 to generate a Type 30 and a Type 1 subtree.
\item If he colored $vu$ (resp. $vr$) with 2 or 3, she may color $vr$ (resp. $vu$) with 3 or 2 to generate a Type 9 and a Type 1 subtree.
\item Suppose $d(r)=3$. If he colored $rx$ with 2 or 3, she may color $vu$ with 2 to generate a Type 28 and a Type 1 subtree, given that $P_4^4$ is forbidden. If he colored an edge of an uncolored $x$-branch, she may color $rx$ with 3 to generate a Type 29 and a Type 1 subtree.
\item If he colored an edge of an uncolored $u$-branch, she may color $vu$ to generate a Type 7 and a Type 2 subtree, given that $S_3^4$ and $P_4^4$ are forbidden.

\end{itemize} 

\noindent\tb{Type 22 subtrees:}

Note that $rs$ and $ry$ are identical by symmetry.

(i) Alice may color $ur$ with 3 to generate a Type 3 and a Type 1 subtree.

(ii) If Bob has colored an edge $e$, Alice may respond as follows.  When $e$ is an edge of an uncolored $r$-branch, our proposed responses of Alice are exactly the same as those for Type 17 subtrees; however, we remark that the permitted types which will be generated by Alice are different.
\begin{itemize}
\item If he colored $wv$ (resp. $ur$) with 4, she may color $ur$ (resp. $wv$) with 4 to generate a Type 1 and a Type 3 subtree.
\item If he colored $vu$ with 4, he generated a Type 3 and a Type 2 subtree.
\item Suppose $d(r)=3$. If he colored $rs$, he generated a Type 3 and a Type 1 subtree. If he colored an edge in an uncolored $s$-branch, she may color $rs$ with an available color in  $\{2,3\}$ to generate a Type 3 and a Type 1 subtree.
\item Suppose $d(r)=4$. If he colored $rs$ with 2 (resp. 3 and 4), she may color $ry$ with 3 (resp. 2 and 2) to generate a Type 3 and a Type 1 subtree.
\item Suppose $d(r)=4$. If he colored an edge in an uncolored $s$-branch. Then, a subtree of Type 22 was split into a Type 1 subtree and the other subtree $H$. If Alice now colors $ur$ with 3 in $H$, she will generate a Type 3 subtree and a 3-LCT $T_r$ with star-node $r$. $T_r$ is of Type 2 if and only if it is not the subtree in Figure~\ref{fig:312}. Therefore, Alice may color $ur$ with 3 if and only if it will not lead to the subtree in Figure~\ref{fig:312}. Note that $T_r$ would be the subtree in Figure~\ref{fig:312} if and only if $d(s)=d(y)=4$ and $sz$ was colored with 2 or 4. So, if $sz$ was colored with 2 (resp. 4), an alternative move of  Alice is to color $ry$ (resp. $rs$) with 2 to generate a Type 25 (resp. Type 26) and a Type 1 subtree.
\end{itemize}

\noindent\tb{Types 23 and 24 subtrees:}

We will provide a unified strategy of Alice for a subtree $H$ of Type 23 or 24.

(i) Alice may color $ur$ with 3 to generate a Type 4 and a Type 2 subtree.

(ii) If Bob has colored an edge, Alice may respond as follows.
\begin{itemize}
\item If he colored an edge of an uncolored $t$-branch, she may color $tw$ with an available color in $\{1,3\}$ to generate a Type 1 and a Type 25 (for a Type 23 $H$) or Type 26 (for a Type 24 $H$) subtree.
\item If he colored $tw$ with 1/3/4, she may color $vu$ with 4 to generate a Type 3 and a Type 9 (for a Type 23 $H$) or Type 27 (for a Type 24 $H$) subtree. Symmetrically, if he colored $vu$ with 4, she may color $tw$ with 1 to generate a Type 1 and a Type 3 subtree.
\item If he colored $wv$ (resp. $ur$) with 4, she may color $ur$ (resp. $wv$) with 4 to generate a Type 3 subtree and a Type 2 (resp. Type 1) subtree.
\item If he colored $rs$ with 3 or 4, he generated a Type 4 and a Type 1 subtree.
\item If he colored an edge of an uncolored $s$-branch, she may color $rs$ with an available color in $\{3,4\}$ to generate a Type 4 subtree and a 2-LCT (for a Type 23 $H$) or 3-LCT (for a Type 24 $H$) with star-node $s$, denoted by $T_s$. $T_s$ is of Type 1 if it has only two colored edges. Suppose $T_s$ has three colored edges. Because $S_3^4$ is forbidden and $r$ was a 4-vertex before $rs$ is colored by Alice, in $T_s$, the 3-SN $s$ is adjacent to at most one 4-vertex, which means $T_s$ is of Type 2.
\end{itemize} 

\noindent\tb{Type 25 and 26 subtrees:}

We will provide a unified strategy of Alice for a subtree $H$ of Type 23 or 24. This unified strategy is very similar to that for Type 23 and 24 subtrees, since the only difference between a Type 23 (resp. Type 24) and a Type 25 (resp. Type 26) subtree is the number of colored edges incident with $w$.

(i) Alice may color $ur$ with 3 to generate a Type 3 and a Type 2 subtree.

(ii) If Bob has colored an edge, Alice may respond as follows.
\begin{itemize}
\item If he colored $wv$ (resp. $ur$) with 4, she may color $ur$ (resp. $wv$) with 4 to generate a Type 3 subtree and a Type 2 (resp. 0) subtree.
\item If he colored $vu$ with 4, he generated a Type 3 and a Type 9 (for a Type 25 $H$) or Type 27 (for a Type 26 $H$) subtree.
\item If he colored $rs$ with 3 or 4, he generated a Type 3 and a Type 1 subtree.
\item If he colored an edge of an uncolored $s$-branch, she may color $rs$ with an available color in $\{3,4\}$  to generate a Type 3 subtree and a 2-LCT (for a Type 25 $H$) or 3-LCT (for a Type 26 $H$) with star-node $s$, denoted by $T_s$. Given that $S_3^4$ is not allowed and using the same argument in the unified strategy for Type 23 and 24 subtrees, $T_s$ is of Type 1 or Type 2.
\end{itemize} 

\noindent\tb{Type 27 subtrees:}

(i) Alice may color $vu$ to generate a Type 3 and a Type 1 subtree.

(ii) If Bob colored $wv$ (resp. an edge of an uncolored $u$-branch), he generated a Type 0 and a Type 2 (resp. a Type 9 and a Type 1) subtree.\\

\noindent\tb{Type 28 subtrees:}

Note that vertices $t$ and $s$ are identical by symmetry.

(i) Alice may color $vu$ with 3 to generate a Type 3 and a Type 1 subtree.

(ii) If Bob has colored an edge, Alice may respond as follows.
\begin{itemize}
\item If he colored $wv$ (resp. $vu$) with 3 or 4, she may color $vu$ (resp. $wv$) with 4 or 3 to generate a Type 3 and a Type 1 subtree.
\item If he colored $vr$ with 3 or 4, she may color $vu$ with 5 to generate a Type 3 and a Type 1 subtree.
\item If he colored $us$ with 3, she may color $wv$ with 3 to generate a Type 0 and a Type 19 subtree.
\item If he colored $us$ with 4 or an edge of the uncolored $s$-branch, she may color $vu$ with 3 to generate a Type 3 and a Type 2 (since $d(s)\leq 3$) subtree.
\end{itemize} 

\noindent\tb{Type 29 subtrees:}

Note that vertices $u$ and $s$ (resp. $w$ and $r$) are identical by symmetry.

(i) Alice may color $vu$ with 2 or 3 to generate a Type 28 (since $P_4^4$ is forbidden) and a Type 1 subtree.

(ii) If Bob colored $wv$ with 3 (resp. 4), Alice may color $vu$ with 4 (resp. 3) to generate a Type 9 and a Type 1 subtree. Symmetrically, If he colored $vu$ with 4, she may color $wv$ with 3 to generate a Type 0 and a Type 9 subtree.

If he colored an edge of an uncolored $u$-branch, she may color $vu$ with an available color in  $\{2,3\}$ to generate a Type 28 and a Type 2 subtree, given that $P_4^4$ is forbidden.\\

\noindent\tb{Type 30 subtrees:}

Note that $d(r),d(s),d(t)\leq 3$ because $S_3^4$ and $P_4^4$ are forbidden.

(i) Alice may color $vu$ with 3 to generate a Type 5 (when $d(r)=2$) or Type 6 (when $d(r)=3$) and a Type 1 subtree.

(ii) If Bob has colored an edge, Alice may respond as follows.
\begin{itemize}
\item If he colored an edge of an uncolored $t$-branch, she may color $tw$ with 3 to generate a Type 1 and a Type 39 subtree.
\item If he colored $tw$ (resp. $vr$) with 3, she may color $vr$ (resp. $tw$) with 3 to generate a Type 40 and a Type 1 subtree.
\item If he colored $wv$ with 3, she may color $us$ with 3 to generate a Type 1 and a Type 5 (if $d(r)=2$) or Type 6 (if $d(r)=3$) subtree.
\item Suppose $d(r)=3$. If he colored $rx$, he generated a Type 1 and a Type 31 subtree. If he colored an edge of an uncolored $x$-branch, she may color $rx$ with 1 or 3 to generate a Type 1 and a Type 31 subtree.
\item If he colored an edge of an uncolored $u$-branch, she may color $vu$ with 3 to generate a Type 2 (since $d(s)\leq 3$) and a Type 5 (if $d(r)=2$) or Type 6 (if $d(r)=3$) subtree.
\end{itemize}

\noindent\tb{Type 31 subtrees:}

Note that $d(t),d(s)\leq 3$, because $P_4^4$ is forbidden.

(i) Alice may color $wv$ with 3 to generate a Type 2 and a Type 9 subtree.

(ii) If Bob has colored an edge, Alice may respond as follows.
\begin{itemize}
\item If he colored an edge of an uncolored $t$-branch, she may color $tw$ with 4 to generate a Type 1 and a Type 32 subtree.
\item If he colored $tw$ with 3 (resp. 4), she may color $vu$ (resp. $wv$) with 3 to generate a Type 3 and a Type 1 (resp. a Type 0 and a Type 9) subtree.
\item If he colored $wv$ (resp. $vu$) with 3 or 4, she may color $vu$ (resp. $wv$) with 4 or 3 to generate a Type 1 (resp. Type 2) and a Type 3 subtree.
\item If he colored $vr$ with color $b=3$ or 4, she may color $tw$ with $b$ to generate a Type 1 and a Type 40 subtree.
\item If he colored $us$ with 3, she may color $wv$ with 3 to generate a Type 2 and a Type 19 subtree.
\item If he colored $us$ with 1/4 or an edge of an uncolored $s$-branch, she may color $vu$ with 3 to generate a Type 9 (since $d(t)\leq 3$) and a Type 2 (since $d(s)\leq 3$) subtree.
\end{itemize}

\noindent\tb{Type 32 subtrees:}

Note that any neighbor of $u$, except $v$, has degree at most 3, because $P_4^4$ is forbidden.

(i) Alice may color $wv$ with 3 to generate a Type 0 and a Type 9 subtree.

(ii) If Bob has colored an edge, Alice may respond as follows.
\begin{itemize}
\item If he colored $wv$ with 5, she may color $vu$ with 3 to generate a Type 3 and a Type 1 subtree.
\item If he colored $vr$ with 3 (resp. 4 or 5), she may color $vu$ with 4 (resp. 3 or 4) to generate a Type 3 and a Type 1 subtree. Similarly, if he colored $vu$ with 3 (resp. 4 or 5), she may color $vr$ with 4 (resp. 5 or 4) to generate a Type 3 and a Type 0 subtree.
\item If he colored $us$ with 3, she may color $vu$ with 5 (if $a=3$) or 4 (if $a=1$) to generate a Type 1 and a Type 33 (if $a=3$) or Type 3 (if $a=1$) subtree.
\item If he colored $us$ with 1 or 4, she may color $vr$ with 4 to generate a Type 0 and a Type 20 subtree.
\item If he colored an edge of an uncolored $s$-branch, she may color $vu$ with 5 (if $a=3$) or 4 (if $a=1$) to generate a Type 2 and a Type 33 (if $a=3$) or 3 (if $a=1$) subtree.
\end{itemize}

\noindent\tb{Type 33 subtrees:}

Coloring anyone of the two uncolored edges will generate a Type 0 and a Type 3 subtree.\\

\noindent\tb{Type 34 subtrees:}

Note that $d(r),d(s)\leq 3$, because $S_3^4$ and $P_4^4$ are forbidden.

(i) Alice may color $tw$ with 3 to generate a Type 1 and a Type 39 subtree.

(ii) If Bob has colored an edge, Alice may respond as follows.
\begin{itemize}
\item If he colored an edge of the uncolored $t$-branch, she may color $tw$ with 3 to generate a Type 1 and a Type 39 subtree.
\item If he colored $wv$ with 3, she may color $us$ with 3 to generate a Type 1 and a Type 5 (if $d(r)=2$) or Type 6 (if $d(r)=3$) subtree.
\item If he colored $vr/vu$ with 3, she may color $tw$ with 3 to generate a Type 1 and a Type 40 subtree.
\item If he colored $rx$ with 1 or 3 (resp. $us$ with 1 or 3), she may color $wv$ with 3 to generate a Type 2 and a Type 9 (resp. a Type 2 and a Type 5 (when $d(r)=2$) or Type 6 (when $d(r)=3$)) subtree.
\item If he colored an edge of an uncolored $x$-branch, she may color $rx$ with an available color in  $\{1,3\}$ to generate a Type 1 and a Type 35 subtree.
\item If he colored an edge of an uncolored $s$-branch, she may color $vu$ with 3 to generate a Type 2 (since $d(s)\leq 3$) and a Type 36 subtree.
\end{itemize}

\noindent\tb{Type 35 subtrees:}

Note that $d(s),d(t)\leq 3$, because $P_4^4$ is forbidden.

(i) Alice may color $wv$ with 3 to generate a Type 2 and a Type 9 subtree.

(ii) If Bob has colored an edge, Alice may respond as follows.
\begin{itemize}
\item If he colored an edge $e$ of the uncolored $t$-branch (only when $d(t)=3$), she may color $wv$ with 3 to generate a Type 3 (if $e$ is incident with $t$) or Type 13 (if $e$ is not incident with $t$) and a Type 9 subtree.
\item If he colored $tw$ with 3, she may color $vu$ with 3 to generate a Type 3 and a Type 1 subtree.
\item If he colored $tw$ with 4, she may color $vr$ with 4 to generate a Type 0 and a Type 40 subtree.
\item If he colored $wv$ (resp. $vu$) with 3 or 4, she may color $vu$ (resp. $wv$) with 4 or 3 to generate a Type 3 subtree, and a Type 1 (resp. 2) subtree.
\item If he colored $vr$ with 3 (resp. 4), she may color $tw$ with 3 (resp. 4) to generate a Type 1 and a Type 40 subtree.
\item If he colored $us$ with 1 or 3, she may color $wv$ with 3 to generate a Type 2 and a Type 19 subtree.
\item If he colored $us$ with 4, she may color $vu$ with 3 to generate a Type 2 and a Type 3 (if $d(t)=2$) or Type 4 (if $d(t)=3$) subtree.
\item If he colored an edge of an uncolored $s$-branch, she may color $vu$ with 3 to generate a Type 2 and a Type 3 (if $d(t)=2$) or Type 4 (if $d(t)=3$) subtree.
\end{itemize}

\noindent\tb{Type 36 subtrees:}

If $d(r)=d(u)=2$, this Type 36 subtree is of Type 3. If $d(r)=2$ and $d(u)=3$, or $d(r)=3$ and $d(u)=2$, this Type 36 subtree is of Type 4. Since $d(r),d(u)\leq 3$, in the following, we assume $d(r)=d(u)=3$.

(i) Alice may color $rw$ with 3 to generate a Type 1 and a Type 4 subtree.

(ii) If Bob has colored an edge, Alice may respond as follows.
\begin{itemize}
\item If he colored an edge of an uncolored $s$-branch (only when $d(s)\geq 2$), she may color $rw$ with 3 to generate a Type 2 and a Type 4 subtree.
\item If he colored $sr$ with 2/3/4, he generated a Type 1 and a Type 4 subtree.
\item If he colored $rw$ (resp. $vu$) with 4, she may color $vu$ (resp. $rw$) with 4 to generate a Type 1 and a Type 3 subtree.
\item If he colored $wv$ with 4, she may color $vu$ with 5 to generate a Type 0 and a Type 1 subtree.
\item If he colored $ut$ with 1/3/4, she may color $wv$ with 4 to generate a Type 2 and a Type 3 subtree.
\item If he colored $us$ with 4, she may color $wv$ with 4 to generate a Type 2 and a Type 3 subtree.
\item If he colored an edge of an uncolored $t$-branch, she may color $ut$ with an available color in $\{1,3\}$ to generate a Type 4 and a Type 1 subtree.
\end{itemize}

\noindent\tb{Type 37 subtrees:}

If $u=r$, this Type 37 subtree is of Type 10. In the following, we may assume $u\neq r$ and let $s$ be the neighbor of $u$ which is not equal to $v$. Since $d(v)=4$ and $P_4^4$ is forbidden, any neighbor of $w$ and $x$ other than $v$ is of degree at most 3. Vertices $w$ and $x$ are identical by symmetry.

(i) Alice may color $vu$ with 2 to generate two Type 2 subtrees.

(ii) If Bob has colored an edge, Alice may respond as follows.
\begin{itemize}
\item If he colored an edge of an uncolored $w$-branch, she may color $wv$ with 2 to generate a Type 2 and a Type 12 subtree.
\item If he colored $wv$, he generated a Type 1 and a Type 12 subtree.
\item If he colored $us$, he generated a Type 10 and a Type 1 subtree.
\item Suppose he colored an edge in an uncolored $s$-branch. Then, this Type 37 subtree was split into a Type 1 subtree and the other subtree, denoted by $H$. If Alice now colors $us$ in $H$, she will generate a Type 10 subtree, and a 2-/3-LCT, denoted by $T_s$. $T_s$ is of Type 2 if it is a 2-LCT. Suppose $T_s$ is a 3-LCT. $T_s$ is of Type 2 if and only if it is not the subtree in Figure~\ref{fig:312} nor~\ref{fig:212}. Therefore, Alice may color $us$ with an appropriate color if and only if it will not lead to these two unpermitted subtrees. In the following, we will propose alternative moves of Alice when coloring $us$ must lead to anyone of these unpermitted subtrees.

\begin{figure}
\centering
\subcaptionbox{Coloring $us$ with 3 would lead to an unpermitted subtree. \label{fig:type37out312}} [.45\linewidth]
 {\includegraphics[scale=.5]{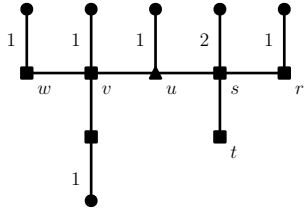}}\hfill%
\subcaptionbox{Coloring $us$ with 3 would lead to an unpermitted subtree. \label{fig:type37out212}} [.45\linewidth]
 {\includegraphics[scale=.5]{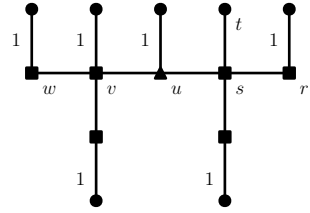}}
	\caption{Situations of Type 37 subtrees in which coloring $us$ with any available color must lead to an unpermitted subtree.}
\end{figure}

If coloring $us$ must generate the subtree in Figure~\ref{fig:312} (resp. Figure~\ref{fig:212}), $H$ can be represented by the configuration in Figure~\ref{fig:type37out312} (resp. Figure~\ref{fig:type37out212}). Alice's alternative move is to color $st$ with 2 to generate a Type 16 (resp. Type 18) and a Type 1 subtree.
\end{itemize}

\noindent\tb{Type 38 subtrees:}

(i) Alice may color $vu$ with 3 to generate two Type 2 subtrees (since $d(t)\leq 3$ when $d(u)=4$).

(ii) If Bob has colored an edge, Alice may respond as follows:
\bit
\item If he colored an edge of the $v$-branch not containing $u$, she may color $vu$ to generate a Type 2 subtree, and a Type 0/3/13 subtree.
\item If he colored $us$ with 1 or 3, she may color $wv$ with 3 to generate a Type 1 and a Type 40 subtree.

%
\item If he colored an edge of an uncolored $s$-branch, she may color $vu$ with any available color to generate two Type 2 subtrees, since $d(t),d(s)\leq 3$.
\item Suppose he colored $ut$ with 1 or 3, she may color $vu$ with 4 to generate a Type 2 and a Type 0 (when $d(u)=3$) or Type 2 (when $d(u)=4$) subtree.
\item Suppose he colored an edge in an uncolored $t$-branch. If $d(u)=4$ and so $d(t)\leq 3$, she may color $ut$ with any available color to generate a Type 6 (since $d(s)\leq 3$) and a Type 1/2 subtree (since $d(t)\leq 3$). 

The strategy of Alice is more complicated if $d(u)=3$ because vertex $t$ can be of degree 4. After Bob's coloring, a Type 37 subtree was split into a Type 1 subtree and the other subtree $H$. If Alice now colors $ut$ in $H$, she will generate a Type 4 subtree, and a 2-/3-LCT, denoted by $T_t$. $T_t$ is of Type 2 if it is a 2-LCT. Suppose $T_t$ is a 3-LCT. $T_t$ is of Type 2 if and only if it is not the subtree in Figure~\ref{fig:312} nor~\ref{fig:212}. Therefore, Alice may color $ut$ with an appropriate color if and only if it will not lead to anyone of these two unpermitted subtrees. In the following, we will propose alternative moves of Alice when coloring $us$ must lead to anyone of these unpermitted subtrees.

\begin{figure}
\centering
\subcaptionbox{Under the constraint that $b\neq 2$. Coloring $ut$ with 1 or 3 would lead to an unpermitted subtree. \label{fig:type38out312}} [.45\linewidth]
 {\includegraphics[scale=.5]{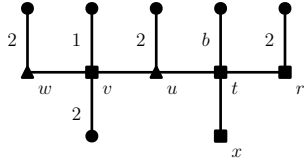}}\hfill%
\subcaptionbox{Coloring $ut$ with 1 or 3 would lead to an unpermitted subtree. \label{fig:type38out212}} [.45\linewidth]
 {\includegraphics[scale=.5]{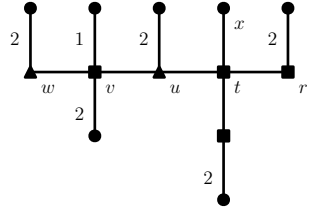}}
	\caption{Situations of Type 38 subtrees in which coloring $ut$ with any available color must lead to an unpermitted subtree.}
\end{figure}

If coloring $ut$ must generate the subtree in Figure~\ref{fig:312} (resp. Figure~\ref{fig:212}), $H$ can be represented by the configuration in Figure~\ref{fig:type38out312} (resp. Figure~\ref{fig:type38out212}). Alice's alternative move is to color $tx$ with 2 to generate a Type 15 (resp. Type 16) and a Type 1 subtree.
\eit
We have presented the strategies of Alice for all Types 1 to 40 in the last and this section so we have completed the proof of Lemma~\ref{lem:inductive}. \eop

\section{Open problem and conjecture}

Since it could be very challenging to prove or disprove that all trees $T$ with $\dt=4$ have $\chi'_g(T)\leq 5$, one may first consider the trees whose 4-vertices form a forest of some basic subclass of trees. For example, one may try to prove the following conjecture proposed by Fong et al.~\citep{FCN}:\\

\noindent\tb{Conjecture} Let $T$ be a finite tree with $\dt=4$. If the subgraph induced by all the 4-vertices is a forest of paths, then $\chi'_g(T)\leq 5$.\\

Moreover, since we confirm in this paper that the 4-vertices are allowed to form a path of three vertices, one may start to consider other subclasses such as stars with three or four leaves:\\

\noindent\tb{Open problem} Let $T$ be a finite tree with $\dt=4$. If the subgraph induced by all the 4-vertices of $T$ is a forest of stars, is $\chi'_g(T)\leq 5$?

\bigskip

\noindent {\bf Acknowledgments} The work of W.H. Chan was partially supported by Guangdong Science and Technology Program Grant 2017A050506025.


\bigskip

\begin{thebibliography}{0}

\bibitem{A2003} S.D. Andres, Spieltheoretische Kantenf{\"a}rbungsprobleme auf W{\"a}ldern und verwandte Strukturen (German), Diplomarbeit, Universit{\"a}t zu K{\"o}ln (2003) 49--50; 55--56; 58; 61; 83--95; 97--123.

\bibitem{A2006} S.D. Andres, The game chromatic index of forests of maximum degree $\Delta\geq 5$, Discrete Applied Mathematics 154 (2006) 1317--1323.

\bibitem{AHS} S.D. Andres, W. Hochst{\"a}ttler, C. Schall{\"u}ck, The game chromatic index of wheels, Discrete Applied Mathematics 159 (2011) 1660--1665.

\bibitem{B}  H.L. Bodlaender, On the complexity of some colouring games, International Journal of Foundations of Computer Science 2 (1991) 133--147.

\bibitem{CZ} L. Cai and X. Zhu, Game chromatic index of $k$-degenerate graphs, Journal of Graph Theory 36 (2001) 144--155.

\bibitem{CN} W.H. Chan and G. Nong. The game chromatic index of some trees of maximum degree 4, Discrete Applied Mathematics 170 (2014) 1--6.

\bibitem{EFHK} P. Erd{\"o}s, U. Faigle, W. Hochst{\"a}ttler and W. Kern, Note on the game chromatic index of trees, Theoretical Computer Science 313 (2004) 371--376.

\bibitem{FC} W.L. Fong, W.H. Chan, The edge coloring game on trees with the number of colors greater than the game chromatic index, Journal of Combinatorial Optimization (2019) 1--25.

\bibitem{FCN} W.L. Fong, W.H. Chan, G. Nong, The game chromatic index of some trees with maximum degree four and adjacent degree-four vertices, Journal of Combinatorial Optimization 36 (2018) 1--12.

\bibitem {Z} X. Zhu, The game coloring number of planar graphs, Journal of Combinatorial Theory, Series B 75 (1999) 245--258.

\end{thebibliography}
\end{document}